\documentclass[oneside,reqno,a4paper,12pt]{amsart}
\usepackage{preamble}

\title{Nakayama functors on proper abelian subcategories}
\author{David Nkansah}
\date{\today}

\begin{document}

\begin{abstract}
    We construct Nakayama functors on proper abelian subcategories of triangulated categories with a Serre functor using approximation theory. This, in turn, allows for the construction of Auslander-Reiten translates. As a result, we prove that suitable proper abelian subcategories are dualising \(k\)-varieties and have enough projectives if and only if they have enough injectives. As an application, we provide a new proof of the existence of Auslander-Reiten sequences in the category of finite dimensional modules over a finite dimensional algebra.
\end{abstract}

\maketitle


{
\footnotesize
\singlespacing
\hypersetup{linkcolor=black}
\tableofcontents
}

\section{Introduction}

    Classical Auslander-Reiten theory was first introduced in \cite{auslander-representation-3-1975} within the context of abelian categories. It was later extended into the triangulated realm in \cite{happel-triangulated-1988}. Relative versions of Auslander-Reiten theory were also introduced in both the abelian \cite{auslander-almost-1981} and triangulated \cite{jorgensen-auslander-2009} cases. This paper aims to study Auslander-Reiten theory in a specific class of abelian subcategories of triangulated categories introduced in \cite{jorgensen-proper-2021}. The class known as the \textit{proper abelian subcategories of triangulated categories} (see Definition~\ref{DEF: proper abelian subcategoires}) generalises hearts of t-structures, and provides access to a theory concerning abelian subcategories, with possible nonzero negative extensions, of negative cluster categories. 
    
    Nakayama functors play an important role in the theory. In the abelian case, they provide an equivalence between the category of projective objects and the category of injective objects and allow for the construction of the Auslander-Reiten translates. Therefore, we begin by seeking a construction of Nakayama functors in proper abelian subcategories. Our results support a change of perspective from hearts of t-structures to proper abelian subcategories, as mentioned in \cite{jorgensen-proper-2021}.

    In this subsection, \(k\) is a field and \(\clT\) is a Krull-Schmidt Hom-finite \(k\)-linear triangulated category with suspension functor \(\Sigma\).
    
\begin{defn}[{\cite[Definition~2.2 (i)]{jorgensen-abelian-2022}}]
    We call a diagram \(X \xrightarrow{x} Y \xrightarrow{y} Z\) in a triangulated category \(\clT\) a \textit{short triangle} if there exists a morphism \(Z \xrightarrow{z} \Sigma X\) such that the augmented diagram \(X \xrightarrow{x} Y \xrightarrow{y} Z \xrightarrow{z} \Sigma X\) is a triangle in \(\clT\).
\end{defn}

\begin{notation}
    We say an exact sequence \(X \xrightarrow{x} Y \xrightarrow{y} Z\) in an abelian category is left (right) exact if \(x\) is a monomorphism (\(y\) is an epimorphism). We say the sequence is short exact if it is both left exact and right exact. In particular, a projective presentation of an object \(X\) in an abelian category is a right exact sequence \(P_1 \xrightarrow{} P_0 \xrightarrow{} X\), where \(P_0\) and \(P_1\) are projective objects in said abelian category. We use similar notation for injective copresentations.
\end{notation}

\begin{defn}
    An \textit{additive subcategory} of an additive category is a full subcategory closed under isomorphisms, direct sums and direct summands.
\end{defn}

\begin{defn}[{\cite[Definition~2.2 (ii)]{jorgensen-abelian-2022}}]\label{DEF: proper abelian subcategoires}
    Let \(\clA\) be an additive subcategory of a triangulated category \(\clT\). Then \(\clA\) is a \textit{proper abelian subcategory of \(\clT\)} if it is abelian and if the following holds:
    \begin{itemize}
        \item The diagram \(X \xrightarrow{x} Y\xrightarrow{y} Z\) is a short exact sequence in \(\clA\) if and only if the diagram \(X \xrightarrow{x} Y\xrightarrow{y} Z\) is a short triangle in \(\clT\) whose terms \(X\), \(Y\) and \(Z\) all lie in \(\clA\).
    \end{itemize}
\end{defn}

\begin{rem}
    Other generalisations of hearts of t-structures include \textit{admissible abelian subcategories} \cite[Definition~1.2.5]{beilinson-faisceaux-1982} and \textit{distinguished abelian subcategories} \cite[Definition~1.1]{linckelmann-abelian-2024}.
\end{rem}

    In the classical case, one may break down Auslander-Reiten theory into the following:
\begin{itemize}
    \item [1.] The study of the collection of non-projective objects and non-injective objects. These are intimately linked to the Auslander-Reiten translates and Auslander-Reiten sequences.
    \item [2.] The study of the collection of projective objects and injective objects. These are intimately linked to Nakayama functors.
\end{itemize}
    Using \cite[Theorem~3.1]{jorgensen-auslander-2009} (see \cite[Theorem~B]{fedele-auslander-2019} for the \((d+2)\)-angulated case), \textit{part 1} is well-understood for proper abelian subcategories: The Auslander-Reiten sequences in a proper abelian subcategory are controlled by approximation properties of the proper abelian subcategory with its ambient triangulated category. However, the aspects mentioned in \textit{part 2} are not covered in \cite{jorgensen-auslander-2009} and therefore, we aim to explore them in this paper. Other considerations about \textit{part 1} can be found in \cite{kleiner-approximations-1997} and \cite{fedele-d-auslander-2020} whereas other considerations about \textit{part 2} can be found in \cite{fedele-almost-2022}.

\begin{defn}[{\cite[page~519]{bondal-representable-1989}}]
    A \textit{Serre functor \(\clX \xrightarrow{\bbS} \clX\)} on a \(k\)-linear category \(\clX\) is a \(k\)-linear autoequivalence together with isomorphisms \(\clX(X,Y) \xrightarrow{\eta_{X,Y}} \D\clX(Y,{\bbS}X)\) that are natural in objects \(X\) and \(Y\) in \(\clX\). Here \((\mod_{k})^{\op} \xrightarrow{\D} \mod_{k}\) denotes the standard \(k\)-dual functor \(\Hom_k(-,k)\) where \(\mod_k\) is the category of finite dimensional right modules over \(k\). In this case, we say \(\clX\) has \textit{Serre duality}.
\end{defn}

    Let \(\clA\) be either the category of finite dimensional (right) modules over a finite dimensional \(k\)-algebra of finite global dimension or a hereditary Ext-finite \(k\)-linear abelian category. Consider \(\clA\) sitting inside its bounded derived category \(\clD^b(\clA)\) as the heart of the canonical t-structure. It was shown in \cite[I.4.6 and Theorem on page~37]{happel-triangulated-1988} and \cite[Corollary~I.3.4]{reiten-noetherian-2002} that a Serre functor on \(\clD^b(\clA)\) restricts to an equivalence between projective objects in \(\clA\) and injective objects in \(\clA\) (see also \cite{happel-auslander-1991}). We would like to point out two problems that may occur in a more general setting. Firstly, the Serre functor might not send a projective object in \(\clA\) to an object in \(\clA\), let alone an injective one (see Example~\ref{EX: A2 into Db(A3)}). Secondly, the Serre functor might not induce an equivalence. Under mild assumptions, Theorem~\ref{Main THM: A} provides a remedy to these issues by constructing a functor via approximations of a Serre functor.

\begin{ex}[Motivational example]\label{EX: A2 into Db(A3)}
    The Auslander-Reiten sequence \(0 \xrightarrow[]{} \substack{1} \xrightarrow[]{} \substack{2\\1} \xrightarrow[]{} \substack{2} \xrightarrow[]{} 0\) in \(\mod{kA_3}\), where \(\substack{2}\) is the unique indecomposable non-projective non-injective \({kA_3}\)-module, induces an exact fully faithful embedding of module categories \(\mod_{kA_2} \xrightarrow{} \mod_{kA_3}\). Therefore, considering \(\mod_{kA_3}\) as the heart of the canonical t-structure on the bounded derived category \(\clD^b(\mod_{kA_3})\), we can identify \(\mod_{kA_2}\) as a full subcategory of \(\clD^b(\mod_{kA_3})\). As the path algebra of the quiver \(A_3\) is of finite global dimension, \(\clD^b(\mod_{kA_3})\) has a Serre functor \({\bbS}\) given by the left derived functor of \((-)\otimes_{kA_3}\D(kA_3)\) (see \cite[Theorem on page 37]{happel-triangulated-1988}). One can check that \({\bbS}\) sends the image of the simple projective \(kA_2\)-module in \(\clD^b(\mod_{kA_3})\) to an object that does not lie in \(\mod_{kA_2}\) (identified as a full subcategory of \(\clD^b(\mod_{kA_3})\)). It is sent to the image of the indecomposable projective-injective object in \(\mod_{kA_3}\). If we consider the Auslander-Reiten quiver of \(\clD^b(\mod_{kA_3})\), the above translates to the mapping \(\substack{1}\xmapsto{{\bbS}}\substack{3\\2\\1}\) in the following diagram:
    
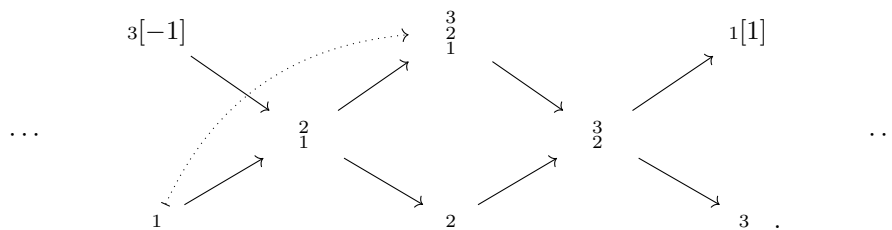
\[
    \begin{tikzcd}
       & \eqmathbox{\substack{3}[-1]} \arrow[rd]                                    &                                       & \eqmathbox{\substack{3\\2\\1}} \arrow[rd] &                                & \eqmathbox{\substack{1}[1]} &        \\
        \cdots &                                                                  & \eqmathbox{\substack{2\\1}} \arrow[rd] \arrow[ru] &                               & \eqmathbox{\substack{3\\2}} \arrow[rd] \arrow[ru] &                   & \cdots \\
       & \eqmathbox{\substack{1}} \arrow[ru] \arrow[rruu, dotted, maps to, bend left] &                                       & \eqmathbox{\substack{2}} \arrow[ru]       &                                       & \eqmathbox{\substack{3}}.      &       
    \end{tikzcd}
    \]
    Although the object \(\S(\substack{1})=\substack{3\\2\\1}\) does not lie in \(\mod_{kA_2}\), we see that there is an irreducible morphism \(\substack{2\\1} \xrightarrow{f} \substack{3\\2\\1}\), where the object \(\substack{2\\1}\) does indeed lie in \(\mod_{kA_2}\). Notice that \(\nu(1) = \substack{2\\1}\), where \(\nu\) is the Nakayama functor on \(\mod_{kA_2}\). We will see in Section~\ref{section: recovering the usual Nakayama functor on modA}, that this morphism \(\nu(\substack{1}) \xrightarrow{f} \S(\substack{1})\) exhibits a strong approximation (see Definition~\ref{defn: (strong) (pre)covers and (pre)envelopes}) of \(\S(\substack{1})\) inside of \(\mod_{kA_2}\).
\end{ex}

    We follow terminology due to \cite[Section 1]{enochs-injective-1981} whereas the alternative terminology in the following is due to Auslander and Smalø, in which the concept was introduced in \cite{auslander-preprojective-1980}.

\begin{defn}\label{defn: (strong) (pre)covers and (pre)envelopes}
    Let \(\clY\) be an additive subcategory of an additive category \(\clX\) and let \(X\) be an object in \(\clX\). For a morphism \(Y \xrightarrow{y} X\) with \(Y\) in \(\clY\), we have the following definitions:
    \begin{itemize}
        \item \(y\) is a \textit{\(\clY\)-precover of \(X\)} if for each object \(Y'\) in \(\clY\) the map \(\clX(Y',Y) \xrightarrow{\clX(Y',y)} \clX(Y',X)\) is an epimorphism. Diagrammatically, this means for each morphism \(Y'\xrightarrow{} X\), there exists a morphism \(Y' \xrightarrow{} Y\) making the following diagram commute:
        \[
        \begin{tikzcd}
                 & Y' \arrow[d] \arrow[ld, dashed, "\exists"'] \\
        Y \arrow[r, "y"] & X .                             
        \end{tikzcd}
        \]
        Alternative terminology for \(y\) is a \textit{right \(\clY\)-approximation of X}.
        \item \(y\) is \textit{right minimal} if each endomorphism \(Y \xrightarrow{y'} Y\) satisfying the equation \(y=yy'\) is automatically an automorphism.
        \item \(y\) is a \textit{\(\clY\)-cover of \(X\)} if it is both a \(\clY\)-precover of \(X\) and right minimal.
        \item \(y\) is a \textit{strong \(\clY\)-cover of \(X\)} if for each object \(Y'\) in \(\clY\) the map \(\clX(Y',Y) \xrightarrow{\clX(Y',y)} \clX(Y',X)\) is an isomorphism. Diagrammatically, this means for each morphism \(Y'\xrightarrow{} X\), there exists a unique morphism \(Y' \xrightarrow{} Y\) making the following diagram commute:
        \[
        \begin{tikzcd}
                 & Y' \arrow[d] \arrow[ld, dashed, "\exists!"'] \\
        Y \arrow[r, "y"] & X.                              
        \end{tikzcd}
        \]
    \end{itemize}
    \textit{\(\clY\)-preenvelopes} (or alternatively, \textit{left \(\clY\)-approximations}), \textit{left minimality, \(\clY\)-envelopes} and \textit{strong \(\clY\)-envelopes} are defined dually.
\end{defn}

\begin{defn}[{\cite[page~307]{auslander-stable-1974}}]\label{DEF: dualising k variety}
    Let \(\clA\) be an essentially small \(k\)-linear category and let \(\Mod_{\clA}\) denote the (abelian) category of \(k\)-linear functors \(\clA^{\op} \xrightarrow{} \Mod_k\) (see \cite[page~184]{auslander-representation-1-1974}). We say that an object \(F\) in \(\Mod_\clA\) is \textit{finitely presented} (see \cite[page~204]{auslander-representation-2-1974}) if there exists a right exact sequence \(\clA(-,Y) \xrightarrow{} \clA(-,X) \xrightarrow{} F\), with \(X\) and \(Y\) objects in \(\clA\) (alternative terminology for such an \(F\) is \textit{coherent}, see \cite[page~189]{auslander-coherent-1966}). Let \(\mod_{\clA}\) denote the full subcategory of \(\Mod_{\clA}\) consisting of the finitely presented objects in \(\Mod_{\clA}\). The standard \(k\)-dual functors induce exact functors
    \begin{equation}\label{eqn: duals}
        (\Mod_{\clA})^{\op} \xrightarrow{\bbD} \Mod_{\clA^{\op}} \text{\qquad and \qquad} \Mod_{\clA^{\op}} \xrightarrow{\bbD} (\Mod_{\clA})^{\op}.
    \end{equation}
    We say \(\clA\) is a \textit{dualising \(k\)-variety} if the functors in (\ref{eqn: duals}) restrict to functors
    \begin{equation*}
        (\mod_{\clA})^{\op} \xrightarrow{\bbD} \mod_{\clA^{\op}} \text{\qquad and \qquad} \mod_{\clA^{\op}} \xrightarrow{\bbD} (\mod_{\clA})^{\op}.
    \end{equation*}
\end{defn}

\subsection{Main results}

    In this subsection, \(k\) is a field and \(\clT\) is a Krull-Schmidt Hom-finite \(k\)-linear triangulated category \(\clT\) with suspension functor \(\Sigma\) and with a Serre functor \(\clT \xrightarrow{\bbS} \clT\). We also fix a full abelian subcategory \(\clA\) of \(\clT\). We let \(\Proj\clA\) denote the full subcategory of \(\clA\) consisting of the projective objects and let \(\Inj\clA\) denote the full subcategory of \(\clA\) consisting of the injective objects.
    
    As witnessed in Example~\ref{EX: A2 into Db(A3)}, the assignment \(P \mapsto {\bbS}P\) for \(P\) a projective object in \(\clA\) only defines a functor \(\Proj \clA \xrightarrow{} \clT\). The idea of the next result is to correct this unwanted feature by requiring the existence of an \(\clA\)-cover of \({\bbS}P\). It turns out, under mild assumptions (see Theorem~\ref{THM: Nakayama functors defined and are equivalances}), this requirement ensures that the assignment \(P \mapsto {\bbS}P \mapsto \nu P\), where \(\nu P \xrightarrow{} {\bbS}P\) is an \(\clA\)-cover of \({\bbS}P\), is functorial and provides an equivalence \(\Proj\clA \simeq \Inj\clA\) of additive categories.

\begin{setup}\label{setup: main theorems}
    Let \(\clA\) be an extension closed \(k\)-linear proper abelian subcategory of \(\clT\) and assume \(\clT(\clA,\Sigma^{-1}\clA)=0\). Further, assume the following:
        \begin{itemize}
        \item For each projective object \(P\) in \(\clA\) there is an \(\clA\)-cover in \(\clT\) of the form \(\nu P\xrightarrow{\alpha_P} {\bbS}P\).
        \item For each injective object \(I\) in \(\clA\) there is an \(\clA\)-envelope in \(\clT\) of the form \({\bbS}^{-1}I\xrightarrow{\beta_I} \nu^{-}I\).
    \end{itemize}
\end{setup}

\begin{ThmIntro}[{Theorem~\ref{THM: Nakayama functors defined and are equivalances}}]\label{Main THM: A}
    Consider Setup~\ref{setup: main theorems}. Then the following hold:
    \begin{itemize}
        \item[1.] The assignment \(P \mapsto \nu P\) augments to an additive functor \(\Proj\clA\xrightarrow{\nu}\Inj\clA\) such that
        \[
        \begin{tikzcd}
            \nu P \arrow[d, "\nu p"'] \arrow[r, "\alpha_P"] & {\bbS}P \arrow[d, "{\bbS}p"] \\
            \nu{P'} \arrow[r, "\alpha_{P'}"]                & {\bbS}P'               
        \end{tikzcd}
        \]
        is commutative for each morphism \(P \xrightarrow{p} P'\) in \(\Proj \clA\).
        \item[2.] The assignment \(I \mapsto \nu^{-} I\) augments to an additive functor \(\Inj\clA\xrightarrow{\nu^-}\Proj\clA\) such that
        \[
        \begin{tikzcd}
            {\bbS}^{-1}I \arrow[r, "\beta_I"] \arrow[d, "{\bbS}^{-1}i"'] & \nu^-I \arrow[d, "\nu^-{i}"] \\
            {\bbS}^{-1}I' \arrow[r, "\beta_{I'}"]                      & \nu^-{I'}                   
        \end{tikzcd}
        \]
        is commutative for each morphism \(I \xrightarrow{i} I'\) in \(\Inj \clA\).
    \end{itemize}
    Moreover, the functors \(\nu\) and \(\nu^-\) are mutual quasi-inverses.
\end{ThmIntro}

    Our next main theorem shows that \(\clA\) enjoys a useful duality condition (see Definition~\ref{DEF: dualising k variety}) and if we impose that \(\clA\) is a length category, then the symmetry between the projective and injective objects in \(\clA\) becomes stronger.

\begin{ThmIntro}[{Theorem~\ref{THM: A is a dualising vartiety} \text{ and } Theorem~\ref{THM: enough proj iff enough inj}}]\label{Main THM: C}
    Consider Setup~\ref{setup: main theorems}. Then the following hold:
    \begin{itemize}
        \item[1.] If \(\clA\) is essentially small and has enough injectives and enough projectives, then \(\clA\) is a dualising \(k\)-variety.
        \item[2.] If every object in \(\clA\) is of finite length, then \(\clA\) has enough projectives if and only if \(\clA\) has enough injectives.
    \end{itemize}
\end{ThmIntro}

    Having Theorem~\ref{Main THM: A} at our disposal, we can construct Auslander-Reiten translates \(\tau\) and \(\tau^-\) in \(\clA\) following the classical pedagogy (see the beginning of Section~\ref{SUBSEC: the construction and fundamental properties} and Definition~\ref{DEF: Auslander-Reiten translates}). The fundamental properties of these mappings are given as our next main result.
    
\begin{ThmIntro}[{Proposition~\ref{PROP: Properties of tau.} \text{ and } Proposition~\ref{PROP: properties of inverse tau}}]\label{Main THM: D}
    Consider Setup~\ref{setup: main theorems} and assume each object in \(\clA\) has a projective cover and an injective envelope. Then there are Auslander-Reiten translates \(\tau\) and \(\tau^-\) on \(\clA\) and they satisfy several standard properties.
\end{ThmIntro}

    The next result is the existence of Auslander-Reiten sequences in \(\clA\).

\begin{ThmIntro}[{Theorem~\ref{THM: Existence of AR sequences in subcateories}}]\label{Main THM: E}
    Consider Setup~\ref{setup: main theorems} and assume each object in \(\clA\) has a projective cover and an injective envelope. Then the following hold:
    \begin{itemize}
        \item [1.] For each indecomposable non-projective object \(C\) in \(\clA\), there exists an Auslander-Reiten sequence in \(\clA\) of the form
            \[
            0 \xrightarrow{} \tau C \xrightarrow{} E_C \xrightarrow{} C \xrightarrow{} 0.
            \]
        \item [2.] For each indecomposable non-injective object \(A\) in \(\clA\), there exists an Auslander-Reiten sequence in \(\clA\) of the form
            \[
            0 \xrightarrow{} A \xrightarrow{} F_A \xrightarrow{} \tau^-A \xrightarrow{} 0.
            \]
    \end{itemize}
\end{ThmIntro}

    When applying the methods presented in the main theorems above to the module category, we recover the standard Nakayama functors and Auslander-Reiten translates (see Theorem~\ref{THM: Nakayama functors in modA} and Proposition~\ref{THM: AR sequences in modA}). As an application, we provide a new proof of the existence of Auslander-Reiten sequences in the finite dimensional module category of a finite dimensional algebra (see Theorem~\ref{THM: existence of AR sequences in mod_A for all A}).

    Theorem~\ref{Main THM: A}, Theorem~\ref{Main THM: C} and Theorem~\ref{Main THM: E} are instances of a fascinating and potentially powerful phenomenon that involves relating intrinsic properties of an abelian category with its relationship to an ambient triangulated category (see the discussion in \cite[page~213]{coelho-simoes-functorially-2022}). In the same paper, it was shown that the heart of a bounded t-structure of a saturated \(\Hom\)-finite Krull-Schmidt \(k\)-linear triangulated category is functorially finite in said triangulated if and only if the heart has enough projective and enough injective objects \cite[Corollary~2.8]{coelho-simoes-functorially-2022} (see \cite[Theorem~2.4]{coelho-simoes-functorially-2022} for a more general statement and compare with Theorem~\ref{Main THM: A} and Theorem~\ref{Main THM: C}, part 2).
    
    Another such example of this phenomenon can be seen in \cite[Proposition~2.6]{jorgensen-abelian-2022}, which says the following: Given an extension closed additive subcategory \(\clA\) of a \(\Hom\)-finite Krull-Schmidt \(k\)-linear triangulated category \(\clT\) such that \(\clT(\clA,\Sigma^{-1}\clA)=0\). If each object of \(\clA\) has a \(\Sigma\clA\)-envelope, then there is an exact structure \(\clE\) on \(\clA\) such that the exact category \((\clA,\clE)\) has enough projective objects.

    \begin{rem}
        In \cite[Remark~2.15]{coelho-simoes-functorially-2022}, the authors mention that it would be interesting to investigate whether their result \cite[Corollary~2.8]{coelho-simoes-functorially-2022} holds without the assumption of the triangulated category being saturated. The results in this paper assume the weaker condition that the triangulated category has a Serre functor (see \cite[page~3]{kawamata-derived-2009}).
    \end{rem}

\subsection{Some useful results}

    For the reader's convenience, we record some results that will be used multiple times throughout this paper. The first result is a useful criterion to use to show when a certain precover is actually a cover.

\begin{lem*}[{\cite[Lemma~2.4]{krause-auslander-2000}}]
    Consider a nonzero morphism \(X \xrightarrow{f} Y\) in an additive category and suppose that \(\End(Y)\) is local. Then \(\alpha\) is left minimal.
\end{lem*}

    The next result allows us to characterise the intrinsic property of an object being projective in a proper abelian subcategory in terms of the ambient triangulated category and its suspension functor.

\begin{thm*}[{\cite[Theorem on page~1]{dyer-exact-2005} (see also \cite[Theorem~3.5]{jorgensen-abelian-2022})}]
     Let \(\clT\) be a triangulated category with suspension functor \(\Sigma\) and let \(\clA\) be an additive full subcategory of \(\clT\) that is closed under extensions and \(\clT(\clA,\Sigma^{-1}\clA)=0\). Then an object \(P\) in \(\clA\) is projective in \(\clA\) if and only if \(\clT(P,\Sigma\clA)=0\).
\end{thm*}

    The final result will be used, sometimes in conjunction with Serre duality, to show that certain \(\Hom\) spaces vanish.

\begin{lem*}[Triangulated Wakamatsu’s Lemma {\cite[Lemma~2.1]{jorgensen-auslander-2009}}]
    Let \(\clT\) be a Krull-Schmidt \(\Hom\)-finite \(k\)-linear triangulated category with suspension functor \(\Sigma\) and let \(\clA\) be a full subcategory of \(\clT\) that is closed under extensions and direct summands. Suppose that a morphism \(A \xrightarrow{\alpha} T\) in \(\clT\) is an \(\clA\)-cover of an object \(T\) in \(\clT\) and complete it to a triangle \(A \xrightarrow{\alpha} T \xrightarrow{} Z \xrightarrow{} \Sigma A\). Then \(\clT(\clA, Z)=0\).
\end{lem*}

\subsection{Global setup}

The following are taken throughout the paper:
\begin{itemize}
    \item \(k\) is a field.
    \item The standard \(k\)-dual functor \((\mod_{k})^{\op} \xrightarrow{\Hom_k(-,k)} \mod_{k}\), where \(\mod_k\) is the category of finite dimensional right modules over \(k\), is denoted by \(\D\).
    \item All subcategories are assumed to be full subcategories closed under isomorphisms.
    \item \(k\)-linear categories are categories enriched over the category of \(k\)-vector spaces with finite direct sums.
    \item For an abelian category \(\clA\), the full subcategory of \(\clA\) consisting of the projective objects is denoted by \(\Proj\clA\) and the full subcategory of \(\clA\) consisting of the injective objects is denoted by \(\Inj\clA\).
    \item \(\clT\) is a Krull-Schmidt Hom-finite \(k\)-linear triangulated category with suspension functor \(\Sigma\) and with a Serre functor \(\clT \xrightarrow{{\bbS}} \clT\).
\end{itemize}

\section{Lemmas}

\subsection{Lemmas on additive categories}
    Consider two objects \(X\) and \(X'\) in an additive category. We denote by \(\Rad(X,X')\) the radical morphisms between \(X\) and \(X'\) as seen in \cite[Lemma~6]{kelly-radical-1964}. That is, a morphism \(X\xrightarrow{x} X'\) in said additive category lies in \(\Rad(X,X')\) if for each morphism \(X' \xrightarrow{x'} X\), we have that \(1_X - x'x\) is an automorphism. The following lemma will be used in Proposition~\ref{PROP: Properties of tau.}.
\begin{lem}\label{LEM: isomorphism between radical maps induced by an additive equivalence}
    Let \(\clX\) and \(\clY\) be additive categories and let \(\clX \xrightarrow{F} \clY\) be a fully faithful additive functor. Then \(F\) induces an isomorphism \(\Rad(X,X')\xrightarrow{}\Rad(FX,FX')\) of abelian groups.
\begin{proof}
    Let \(X \xrightarrow{x} X'\) be a radical morphism in \(\clX\). We will first show that \(FX \xrightarrow{Fx} FX'\) is a radical morphism in \(\clY\). To this end, let \(FX' \xrightarrow{g} FX\) be a morphism in \(\clY\). As \(F\) is fully faithful, \(g\) is of the form \(Fx'\) for some morphism \(X' \xrightarrow{x'} X\) in \(\clX\). As \(x\) is a radical morphism in \(\clX\), we have that \(1_X-x'x\) is an automorphism. Therefore, \(1_{FX}-gF(x)\) is automorphism. We now have that the group isomorphism \(\clX(X,X') \xrightarrow{} \clY(FX,FX')\) induced by \(F\) restricts to an injective group homomorphism \(\Rad(X,X') \xrightarrow{} \Rad(FX, FX')\). For surjectivity, we let \(FX \xrightarrow{f} FX'\) be a radical morphism in \(\clY\). Again, we may take \(f\) to be of the form \(Fw\) for some morphism \(X \xrightarrow{w} X'\) in \(\clX\). It suffices to show that \(w\) is a radical morphism in \(\clX\). Let \(X' \xrightarrow{w'} X\) be a morphism in \(\clX\). Then \(F(1_{X}-w'w)=1_{FX}-F(w')f\) is an automorphism as \(f\) is a radical morphism in \(\clY\). As \(F\) is fully faithful, it reflects isomorphisms. Therefore, \(1_{X}-w'w\) is an automorphism and hence \(w\) is a radical morphism in \(\clX\).
\end{proof}
\end{lem}

    The following Lemma will be a diagram trick involving a strong cover. It will be used in Lemma~\ref{REM: any choice of strong covers yield nat isom fucntors} and Theorem~\ref{THM: Nakayama functors defined and are equivalances}.

\begin{lem}\label{lem: back square commutes with strong cover}
    Let \(\clX\) be an additive category, let \(\clY\) be an additive subcategory of \(\clX\) and consider the following diagram in \(\clX\):
        \[
        \begin{tikzcd}
            Y_1 \arrow[rr, "a"] \arrow[rd, "e"'] \arrow[dd, "b"'] &                                 & X_1 \arrow[ld, "f"] \arrow[dd, "c"] \\
                                                                             & X_2 \arrow[dd, "h"  {yshift=16pt}] &                                                  \\
        X_3 \arrow[rr, "d"{xshift=-20pt}] \arrow[rd, "g"']                  &                                 & Y \arrow[ld, "\alpha"]                     \\
                                                                             & X,                              &                                                 
        \end{tikzcd}
    \]
    with the following properties:
    \begin{itemize}
        \item The objects \(Y\) and \(Y_1\) both lie in \(\clY\).
        \item All but the back-most square with vertices \(Y_1, X_1, X_3\) and \(Y\) commute.
    \end{itemize}
    If \(\alpha\) is a strong \(\clY\)-cover of \(X\), then the back-most square commutes.
\begin{proof}
    By the commutativity we have, the following equalities hold \(\alpha d b = gb = he = hfa =\alpha ca\). So if \(\alpha\) is a strong \(\clY\)-cover, then we get \(db = ca\) as required.
\end{proof}
\end{lem}

    The following lemma will be used multiple times throughout this paper, where the functor \(\clX \xrightarrow{S} \clX\) in the statement of the lemma will taken to be the Serre functor \(\clT \xrightarrow{\bbS} \clT\) mentioned in the Global setup. It will allow us to construct additive functors through the existence of approximations.
    
\begin{lem}\label{LEM: Nakayama functor defined 0}
    Let \(\clX\) be an additive category and let \(\clX \xrightarrow{S} \clX\) be an additive endofunctor. Let \(\clY\) be an additive subcategory of \(\clX\) and let \(\clZ\) be an additive subcategory of \(\clY\). Then the following hold:
    \begin{itemize}
        \item[1.] If for each object \(Z\) in \(\clZ\) there is a strong \(\clY\)-cover in \(\clX\) of the form \(C_Z\xrightarrow{c_Z} SZ\), then the assignment \(Z\mapsto C_Z\) on objects augments to an additive functor \(\clZ\xrightarrow{C_{(-)}}\clY\) such that the induced diagram
        \[
        \begin{tikzcd}
            C_Z \arrow[d, "C_{z}"'] \arrow[r, "c_{Z}"] & SZ \arrow[d, "Sz"] \\
            C_{Z'} \arrow[r, "c_{Z'}"]                 & SZ'               
        \end{tikzcd}
        \]
        is commutative for each morphism \(Z \xrightarrow{z} Z'\) in \(\clZ\).
        \item[2.] If for each object \(Z\) in \(\clZ\) there is a strong \(\clY\)-envelope in \(\clX\) of the form \(S^{-1}Z\xrightarrow{} E_Z\), then the assignment \(Z\mapsto E_Z\) on objects augments to an additive functor \(\clZ\xrightarrow{}\clY\) such that the induced diagram
        \[
        \begin{tikzcd}
            S^{-1}Z \arrow[r] \arrow[d, "S^{-1}z"'] & E_Z \arrow[d, "E_{z}"] \\
            S^{-1}Z' \arrow[r]                      & E_{Z'}                
        \end{tikzcd}
        \]
        is commutative for each morphism \(Z \xrightarrow{z} Z'\) in \(\clZ\).
    \end{itemize}
\begin{proof}
    \textit{Part 1.} Given a morphism \(Z \xrightarrow{z} Z'\) between objects in \(\clZ\), we have the composition \(C_Z \xrightarrow{c_Z} SZ \xrightarrow{Sz} SZ'\). As \(C_{Z'} \xrightarrow{c_{Z'}} S{Z'}\) is a strong \(\clY\)-cover in \(\clX\), there exists a unique morphism \(C_Z \xrightarrow{C_z} C_{Z'}\) such that \(S(z)c_Z=c_{Z'}C_z\). The uniqueness of \(C_z\) in the assignment \(z \mapsto C_z\) on morphisms ensures functoriality. 
    
    For additivity, let \(Z \xrightarrow{z_1,z_2} {Z'}\) be morphisms between objects in \(\clZ\). We need to show \(C_{({z_1+z_2})}=C_{z_1}+C_{z_2}\). By their respective definitions and as \(S\) is an additive functor, we have the following equalities:
    \[
        c_{Z'}C_{({z_1+z_2})} = S(z_1+z_2)c_Z = (Sz_1+Sz_2)c_Z = c_{Z'}(C_{z_1} + C_{z_2}).
    \]
    Additivity follows as \(c_{Z'}\) is a strong \(\clY\)-cover.

    \textit{Part 2.} Dual to \textit{part 1}.
\end{proof}
\end{lem}

    Lemma \ref{LEM: Nakayama functor defined 0} tells us that a global choice of strong \(\clY\)-covers of the objects \(SZ\), for each object \(Z\) in \(\clZ\), give rise to an additive functor \(\clZ \xrightarrow{} \clY\) fitting into the commutative diagram in said lemma. The next lemma, however, uses the uniqueness of the \(\clY\)-covers to ensure that any other global choice of \(\clY\)-covers of the objects \(SZ\) would result in a canonically naturally isomorphic functor. It will be used in Theorem~\ref{THM: Nakayama functors in modA}.

\begin{lem}\label{REM: any choice of strong covers yield nat isom fucntors}
    Let \(\clX\) be an additive category and let \(\clX \xrightarrow{S} \clX\) be an additive endofunctor. Let \(\clY\) be an additive subcategory of \(\clX\) and let \(\clZ\) be an additive subcategory of \(\clY\). Then the following hold:
    \begin{itemize}
        \item[1.] Assume for each object \(Z\) in \(\clZ\) there are strong \(\clY\)-covers in \(\clX\) of the form \(C_Z\xrightarrow{c_Z} SZ\) and \(D_Z\xrightarrow{d_Z} SZ\). The assignments \(Z \mapsto C_Z\) and \(Z \mapsto D_Z\) augment to two functors \(\clZ \xrightarrow{} \clY\) by Lemma~\ref{LEM: Nakayama functor defined 0}. These functors are naturally isomorphic.
        \item[2.] Assume for each object \(Z\) in \(\clZ\) there are strong \(\clY\)-envelopes in \(\clX\) of the form \(S^{-1}Z\xrightarrow{} E_Z\) and \(S^{-1}Z\xrightarrow{} F_Z\). The assignments \(Z \mapsto E_Z\) and \(Z \mapsto F_Z\) augment to two functors \(\clZ \xrightarrow{} \clY\) by Lemma~\ref{LEM: Nakayama functor defined 0}. These functors are naturally isomorphic.
    \end{itemize}
\begin{proof}
    Assume for each object \(Z\) in \(\clZ\) we have strong \(\clY\)-covers in \(\clX\) of the form \(C_Z\xrightarrow{c_Z} SZ\) and \(D_Z\xrightarrow{d_Z} SZ\). By Lemma~\ref{LEM: Nakayama functor defined 0}, we have two functors \(\clZ \xrightarrow{} \clY\), given by the assignments \(Z \mapsto C_Z\) and \(Z \mapsto D_Z\). As \(C_Z\) is an object in \(\clY\) and \(d_Z\) is a strong \(\clY\)-cover, there exists a unique isomorphism \(C_Z \xrightarrow{\varphi_Z} D_Z\) such that \(c_Z=d_Z\varphi_Z\). We now show that the collection \(\{\varphi_Z \mid Z \text{ an object in } \clZ\}\) of isomorphisms in \(\clY\) form the components of a natural isomorphism between the functors in question. Let \(Z \xrightarrow{z} Z'\) be a morphism in \(\clZ\). We obtain the diagram:
    \[
\begin{tikzcd}
C_Z \arrow[rr, "\varphi_Z"] \arrow[dd, "C_z"'] \arrow[rd, "c_Z"'] &                     & D_Z \arrow[dd, "D_z"] \arrow[ld, "d_Z"] \\
                                                                  & SZ \arrow[dd, "Sz" {yshift=16pt}] &                                         \\
C_{Z'} \arrow[rr, "\varphi_{Z'}" {xshift=-20pt}] \arrow[rd, "c_{Z'}"']                    &                     & D_{Z'}, \arrow[ld, "d_{Z'}"]                   \\
                                                                  & S{Z'}                  &                                        
\end{tikzcd}
    \]
    where all but the back-most square with vertices \(C_Z, C_{Z'}, D_{Z'}\) and \(D_Z\) commute. By the application of Lemma~\ref{lem: back square commutes with strong cover}, we are done.
\end{proof}
\end{lem}

\subsection{A lemma on \texorpdfstring{\(k\)}{k}-linear abelian subcategories}

    The following lemma will be used in Proposition~\ref{PROP: equivalent definition of nu}.

\begin{lem}\label{LEM: strong covers are representables in an abelian subcat}
    Let \(\clA\) be a \(k\)-linear abelian subcategory of \(\clT\). Let \(X\) and \(P\) be objects in \(\clA\) with \(P\) projective. Then the following two statements are equivalent:
    \begin{itemize}
        \item[1.] There is a strong \(\clA\)-cover in \(\clT\) of the form \(X \xrightarrow{\alpha_P} {\bbS}P\).
        \item[2.] There is a natural isomorphism of the form \(\D\clA(P,-) \cong \clA(-, X)\) of functors \(\clA^{\op} \xrightarrow{} \mod_k\), i.e. the object \(X\) represents the functor \(\D\clA(P,-)\).
    \end{itemize}
    Dually, let \(Y\) and \(I\) be objects in \(\clA\) with \(I\) injective. Then the following two statements are equivalent:
    \begin{itemize}
        \item[1'.] There is a strong \(\clA\)-envelope in \(\clT\) of the form \({\bbS}^{-1}I \xrightarrow{} Y\).
        \item[2'.] There is a natural isomorphism of the form \( \D\clA(-,I) \cong \clA(Y,-)\) of functors \(\clA \xrightarrow{} \mod_k\), i.e. the object \(Y\) represents the functor \(\D\clA(-,I)\).
    \end{itemize}
\begin{proof}
    \((1\Rightarrow 2)\): Suppose there is a strong \(\clA\)-cover of the form \(X \xrightarrow{\alpha_P} {\bbS}P\). We have a natural isomorphism \(\clA(-,X) \xrightarrow{\restr{\clT(-,\alpha_P)}{\clA}} \restr{\clT(-,{\bbS}P)}{\clA} \cong \D\clA(P,-)\) of functors \(\clA^{\op} \xrightarrow{} \mod_k\), where the first whiskered composite natural transformation is a natural isomorphism as \(\alpha_P\) is a strong \(\clA\)-cover and the second natural isomorphism is given by Serre duality. Therefore, \(X\) represents the functor \(\D\clA(P,-)\).

    \((2\Rightarrow 1)\): Suppose the object \(X\) represents the functor \(\D\clA(P,-)\). Then we have a natural isomorphism \(\gamma \colon \clA(-,X) \cong \D\clA(P,-) \cong \restr{\clT(-,{\bbS}P)}{\clA}\) of functors \(\clA^{\op} \xrightarrow{} \mod_k\), where the second natural isomorphism is given by Serre duality. We will show that the morphism \(X \xrightarrow{\gamma_{X}(1_{X})} {\bbS}P\), where \(\gamma_{X}\) is the component of the natural isomorphism \(\gamma\) at \(X\), is a strong \(\clA\)-cover. Let \(\alpha_P\coloneq\gamma_{X}(1_{X})\). It suffices to show for each object \(A\) in \(\clA\), the component \(\gamma_A\) coincides with \(\clA(A,X) \xrightarrow{\restr{\clT(A,\alpha_P)}{\clA}} \restr{\clT(A,{\bbS}P)}{\clA}\). To this end, let \(A \xrightarrow{a} X\) be a morphism in \(\clA\). By naturality, the following diagram is commutative:
    \[
    \begin{tikzcd}
        {\clA(X,X)} \arrow[r, "\gamma_{X}"] \arrow[d, "{\clA(a,X)}"'] & {\restr{\clT(X,{\bbS}P)}{\clA}} \arrow[d, "{\restr{\clT(a,{\bbS}P)}{\clA}}"] \\
        {\clA(A,X)} \arrow[r, "\gamma_A"]                          & {\restr{\clT(A,{\bbS}P)}{\clA}}.                         
    \end{tikzcd}
    \]
    Chasing the identity \(1_{X}\) through this diagram yields the equality \(\gamma_A(a)=\restr{\clT(A,\alpha_P)}{\clA}(a)\).

    \textit{Part 2.} Dual to \textit{part 1}.
\end{proof}
\end{lem}

\subsection{A lemma on abelian length categories}

    The following lemma will be used in Theorem~\ref{THM: enough proj iff enough inj}.

\begin{lem}\label{LEM: finite length categories, simples and enough injectives/projectives}
    Let \(\clA\) be an abelian category for which each object of \(\clA\) is of finite length. Then the following hold:
    \begin{itemize}
        \item[1.] If each simple object in \(\clA\) admits a nonzero morphism to an injective object in \(\clA\), then \(\clA\) has enough injectives.
        \item[2.] If each simple object in \(\clA\) admits a nonzero morphism from a projective object in \(\clA\), then \(\clA\) has enough projectives.
    \end{itemize}
\begin{proof}
    \textit{Part 1.} Let \(A\) be an object in \(\clA\). As \(A\) is of finite length, there exists a finite composition series \(0=A_0 \subseteq A_1 \subseteq \cdots \subseteq A_{n-1} \subseteq A_n=A\) of \(A\). Denote the simple quotients \(A_{i+1}/A_i\) by \(S_{i+1}\) and note \(A_1=S_1\). There are short exact sequences of the form \(A_i \xrightarrow{} A_{i+1} \xrightarrow{} S_{i+1}\) and by assumption, there are nonzero morphisms \(S_{i+1}\xrightarrow{} I_{i+1} \), where every \(I_{i+1}\) is an injective object in \(\clA\). In particular, these morphisms are monomorphisms in \(\clA\). By the dual of the Horseshoe Lemma \cite[Lemma~8.2.1]{enochs-relative-2011}, we inductively get monomorphisms \(A_{i+1} \xrightarrow{} I_1\oplus I_2\oplus \cdots \oplus I_{i+1}\) in \(\clA\). In particular, we acquire a monomorphism \(A=A_n \xrightarrow{} I_1\oplus I_2\oplus \cdots \oplus I_n\) in \(\clA\), proving our statement.

    \textit{Part 2.} Dual to \textit{part 1}.
\end{proof}
\end{lem}

\subsection{A lemma on proper abelian subcategories}

    The following lemma will be used in Proposition~\ref{PROP: Properties of tau.}, Theorem~\ref{THM: Existence of AR sequences in subcateories} and implicitly in Proposition~\ref{Prop: certain covers and envelops in modA}.

\begin{lem}\label{LEM: A specific precover and preenvelope}
    Let \(\clA\) be an extension closed \(k\)-linear proper abelian subcategory of \(\clT\) and assume \(\clT(\clA,\Sigma^{-1}\clA)=0\). Further, assume the following:
    \begin{itemize}
        \item For each projective object \(P\) in \(\clA\) there is a strong \(\clA\)-cover in \(\clT\) of the form \(C_P\xrightarrow{c_P} {\bbS}P\).
        \item For each injective object \(I\) in \(\clA\) there is a strong \(\clA\)-envelope in \(\clT\) of the form \({\bbS}^{-1}I\xrightarrow{} E_I\).
    \end{itemize}
    As \(\Proj \clA\) is an additive subcategory of \(\clA\), Lemma~\ref{LEM: Nakayama functor defined 0} gives rise to a functor \(\Proj\clA \xrightarrow{} \clA\) given by the assignment \(P \mapsto C_P\). The following then hold:
    \begin{itemize}
        \item[1.] If \(P_1 \xrightarrow{p_1} P_0 \xrightarrow{p_0} A\) is a projective presentation of \(A\) in \(\clA\), then there is an \(\clA\)-precover of the form \(X_A \xrightarrow{\theta} \Sigma^{-1}{\bbS}A\), where \(X_A\) is the kernel of the induced morphism \(C_{P_1} \xrightarrow{C_{p_1}} C_{P_0}\) in \(\clA\). Moreover, if \(A\) non-projective object in \(\clA\), then the \(\clA\)-precover \(\theta\) is necessarily nonzero.
        
        \item[2.] If \(B \xrightarrow{} I^0 \xrightarrow{} I^1\) is an injective copresentation of \(B\) in \(\clA\), then there is an \(\clA\)-preenvelope of the form \(\Sigma {\bbS}^{-1}B \xrightarrow{\xi} Z_B\), where \(Z_B\) is the cokernel of the induced morphism \(E_{I^0} \xrightarrow{} E_{I^1}\) in \(\clA\). Moreover, if \(B\) non-injective object in \(\clA\), then the \(\clA\)-preenvelope \(\xi\) is necessarily nonzero.
    \end{itemize}
\begin{proof}
    \textit{Part 1.} Our assumption allows us to use Lemma~\ref{LEM: Nakayama functor defined 0} to produce an additive functor \(\Proj\clA \xrightarrow{} \clA\) given by \(P \mapsto C_P\). The projective presentation \(P_1 \xrightarrow{p_1} P_0 \xrightarrow{p_0} A\) of \(A\) induces a left exact sequence \(X_A \xrightarrow{x_A} C_{P_1} \xrightarrow{C_{p_1}} C_{P_0}\), where \(x_A\) is the kernel of \(C_{p_1}\). Complete the morphism \({\bbS}P_1 \xrightarrow{{\bbS}p_1} {\bbS}P_0\) to a short triangle \(Y \xrightarrow{y} {\bbS}P_1 \xrightarrow{{\bbS}p_1} {\bbS}P_0\). As \(\clA\) is a proper abelian subcategory, the short exact sequence \(\Omega A \xrightarrow{\iota_1} P_0 \xrightarrow{p_0} A\), where \(\Omega A\) is the syzygy of \(A\), is a short triangle. By rotating, \(\Sigma^{-1}P_0 \xrightarrow{-\Sigma p_0} \Sigma^{-1}A \xrightarrow{} \Omega A\) is a short triangle. Therefore, there is a short triangle \(\Sigma^{-1}{\bbS}P_0 \xrightarrow{-\Sigma {\bbS}p_0} \Sigma^{-1}{\bbS}A \xrightarrow{w} {\bbS}\Omega A\). Here we used that \({\bbS}\) is a triangulated functor (see \cite[Proposition~3.3]{bondal-representable-1989} for the classical reference and see \cite[Theorem~A.4.4]{bocklandt-graded-2008} for a proof which makes the natural transformation more evident). After rotating triangles, we get the following solid commutative diagram:
    \begin{equation}\label{D2}
    \begin{tikzcd}[column sep =8ex]
        0 \arrow[r]                        & X_A \arrow[r, "x_A"] \arrow[d, dashed, "\theta_1"]    & C_{P_1} \arrow[r, "C_{p_1}"] \arrow[d, "c_{P_1}"] & C_{P_0} \arrow[d, "c_{P_0}"]\\
        \Sigma^{-1}{\bbS}P_0 \arrow[r] \arrow[d, equal] & Y \arrow[r, "y"] \arrow[d, dashed, "\theta_2"]      & {\bbS}P_1 \arrow[r, "{\bbS}p_1"] \arrow[d, "{\bbS}\pi_1"] & {\bbS}P_0 \arrow[d, equal]\\
        \Sigma^{-1}{\bbS}P_0 \arrow[r, "-\Sigma {\bbS}p_0"]          & \Sigma^{-1}{\bbS}A \arrow[r, "w"] & {\bbS}\Omega A \arrow[r, "{\bbS}\iota_1"] & {\bbS}P_0
    \end{tikzcd}
    \end{equation}
    where the composition \({\bbS}P_1 \xrightarrow{{\bbS}\pi_1} {\bbS}\Omega A \xrightarrow{{\bbS}\iota_1} {\bbS}P_0\) is obtained by applying \({\bbS}\) to the usual factorisation \(p_1\colon P_1 \xrightarrow{\pi_1} \Omega A \xrightarrow{\iota_1} P_0\) of the morphism \(p_1\). Note that the bottom and the middle rows are both triangles. By an axiom of triangulated categories, there exists a morphism \(Y \xrightarrow{\theta_2} \Sigma^{-1}{\bbS}A\). Since \({\bbS}(p_1)c_{P_1}x_A=0\), there exists a morphism \(X_A \xrightarrow{\theta_1} Y\) as \(y\) is a weak kernel of \({\bbS}p_1\). The morphisms \(\theta_1\) and \(\theta_2\) make the entire diagram~\ref{D2} commute. Taking \(\theta=\theta_2\theta_1\), we show that \(\theta\) is indeed an \(\clA\)-precover. The next steps in our proof will consist of adding many new morphisms to diagram~\ref{D2}, the reader is suggested to use diagram~\ref{D3} below to aid their reading (the dashed arrows in diagram~\ref{D3} will be introduced in the following argument).
        \begin{equation}\label{D3}
    \begin{tikzcd}[row sep= normal, column sep =8ex]
        0 \arrow[r]                                     & X_A \arrow[rr, "x_A"] \arrow[dd, "\theta_1"] &                                                                                         & C_{P_1} \arrow[r, "C_{p_1}"] \arrow[dd, "c_{P_1}"]   & C_{P_0} \arrow[dd, "c_{P_0}"] \\
                                                &                                                      & B \arrow[lu, "\overline{b_1}"', dashed] \arrow[ld, "b"', dashed] \arrow[ru, "b_1", dashed]                       &                                                      &                               \\
        \Sigma^{-1}{\bbS}P_0 \arrow[r] \arrow[dd, equal]         & Y \arrow[rr, "y"] \arrow[dd, "\theta_2"]     &                                                                                         & {\bbS}P_1 \arrow[r, "{\bbS}p_1"] \arrow[dd, "{\bbS}\pi_1"] & {\bbS}P_0 \arrow[dd, equal]            \\
                                                &                                                      & B \arrow[lu, "\overline{b_3}"', dashed] \arrow[ld, "b_2"', dashed] \arrow[ru, "b_3", dashed] \arrow[rd, "wb_2", dashed] &                                                      &                               \\
        \Sigma^{-1}{\bbS}P_0 \arrow[r, "-\Sigma {\bbS}p_0"] & \Sigma^{-1}{\bbS}A \arrow[rr, "w"]                     &                                                                                         & {\bbS}\Omega A \arrow[r, "{\bbS}\iota_1"]                & {\bbS}P_0.                     
    \end{tikzcd}
    \end{equation}
    First, we show that \(\theta_1\) is an \(\clA\)-precover. Let \(B \xrightarrow{b} Y\) be a morphism with the object \(B\) in \(\clA\). As \(c_{P_1}\) is a strong \(\clA\)-precover, there exists a unique morphism \(B \xrightarrow{b_1} C_{P_1}\) such that \(yb=c_{P_1}b_1\). We have the following equalities: \(c_{P_0}C_{p_1}b_1={\bbS}(p_1)c_{P_1}b_1={\bbS}(p_1)yb=0\), where the last equality holds as consecutive morphisms in a triangle vanish. Then \(C_{p_1}b_1=0\) as \(c_{P_0}\) is a strong \(\clA\)-cover. As \(x_A\) is the kernel of \(C_{p_1}\), there exists a unique morphism \(B \xrightarrow{\overline{b_1}} X_A\) such that \(b_1=x_A\overline{b_1}\). We have the following equalities: \(yb=c_{P_1}b_1=c_{P_1}x_A\overline{b_1}=y\theta_1\overline{b_1}\), which gives \(y(b-\theta_1\overline{b_1})=0\). As \(y\) is a weak kernel of \({\bbS}{p_1}\), the morphism \(b-\theta_1\overline{b_1}\) factors through a morphism in the \(k\)-vector space \(\clT(B,\Sigma^{-1}{\bbS}P_0)\). But, \(\clT(B,\Sigma^{-1}{\bbS}P_0)\cong\D\clT(P_0,\Sigma B)=0\) where the isomorphism is given by Serre duality and where the equality holds as \(P_0\) is a projective object in \(\clA\) (see \cite[Theorem on page~1]{dyer-exact-2005}). Hence, we are done.

    It now suffices to show that the map \(\clT(B,Y) \xrightarrow{\clT(B,\theta_2)} \clT(B,\Sigma^{-1}{\bbS}A)\) is an epimorphism for each object \(B\) in \(\clA\). To this end, let \(B \xrightarrow{b_2} \Sigma^{-1}{\bbS}A\) be a morphism with \(B\) in \(\clA\). By Serre duality, we have the following commutative diagram:
    \[
    \begin{tikzcd}[row sep = 5ex]
        {\clT(B,{\bbS}P_1)} \arrow[r, "\cong"] \arrow[d, "{\clT(B,{\bbS}\pi_1)}"'] & {\D\clT(P_1,B)} \arrow[d, "{\D\clT(\pi_1,B)}"] \\
        {\clT(B,{\bbS}\Omega A)} \arrow[r, "\cong"]                          & {\D\clT(\Omega A,B)}.                         
    \end{tikzcd}
    \]
    As \(\pi_1\) is an epimorphism, the map \(\D\clT(\pi_1,B)\) is an epimorphism and hence, so is the map \(\clT(B,{\bbS}\pi_1)\). Therefore, there exists a morphism \(B \xrightarrow{b_3} {\bbS}P_1\) such that \(wb_2={\bbS}(\pi_1)b_3\). We have the following equalities: \({\bbS}(p_1)b_3={\bbS}(\iota_1){\bbS}(\pi_1)b_3={\bbS}(\iota_1)wb_2=0\), where the last equality holds as consecutive morphisms in a triangle vanish. As \(y\) is a weak kernel of \({\bbS}p_1\), there exists a morphism \(B \xrightarrow{\overline{b_3}} Y\) such that \(b_3=y\overline{b_3}\). We have the following equalities: \(wb_2={\bbS}(\pi_1)b_3={\bbS}(\pi_1)y\overline{b_3}=w\theta_2\overline{b_3}\), which gives \(w(b_2-\theta_2\overline{b_3})=0\). As \(w\) is a weak kernel of \({\bbS}\iota_1\), the morphism \(b_2-\theta_2\overline{b_3}\) factors through a morphism in the \(k\)-vector space \(\clT(B,\Sigma^{-1}{\bbS}P_0)\). But, \(\clT(B,\Sigma^{-1}{\bbS}P_0)\cong\D\clT(P_0,\Sigma B)=0\) where the isomorphism is given by Serre duality and where the equality holds as \(P_0\) is a projective object in \(\clA\) (see \cite[Theorem on page~1]{dyer-exact-2005}). Hence, we are done.

    Suppose that \(A\) is a non-projective object in \(\clA\). By \cite[Theorem on page~1]{dyer-exact-2005}, there exists an object \(T\) in \(\clA\) with \(\clT(A,\Sigma T)\neq0\). Then the \(k\)-vector space \(\D\clT(A,\Sigma T)\) is also nonzero and we have \(\clT(T,\Sigma^{-1}\bbS A) \cong \D\clT(A,\Sigma T) \neq 0\), where the isomorphism holds by Serre duality. Therefore, we can pick a nonzero morphism \(T \xrightarrow{t} \Sigma^{-1}\bbS A\). As \(\theta\) is an \(\clA\)-precover, the nonzero morphism \(t\) must factors through \(\theta\), resulting in \(\theta\) necessarily being nonzero.
    
    \textit{Part 2.} Dual to \textit{part 1}.
\end{proof}
\end{lem}

\section{Nakayama functors on proper abelian subcategories} \label{SECTION: main theory}

\subsection{The construction and fundamental properties}

\begin{setup}\label{SETUP: Base of proper abelian}
    We fix \(\clA\) to be an extension closed \(k\)-linear proper abelian subcategory of \(\clT\) and assume \(\clT(\clA,\Sigma^{-1}\clA)=0\).
\end{setup}

\begin{prop}\label{PROP: strong covers and envelopes}
    Let \(\clA\) be a category with Setup~\ref{SETUP: Base of proper abelian}. For \(P\) a projective object in \(\clA\) and \(I\) an injective object in \(\clA\), the following hold:
    \begin{itemize}
        \item[1.] Each \(\clA\)-cover in \(\clT\) of the form \(C_P \xrightarrow{ c_P} {\bbS}P\) is a strong \(\clA\)-cover in \(\clT\).
        \item[2.] Each \(\clA\)-envelope in \(\clT\) of the form \({\bbS}^{-1}I \xrightarrow{} E_I\) is a strong \(\clA\)-envelope in \(\clT\).
    \end{itemize}
\begin{proof}
    Let \(A\xrightarrow{a} C_P\) be a morphism in \(\clA\) such that \( c_Pa=0\). In order to show that \( c_P\) is a strong \(\clA\)-cover in \(\clT\) we need \(a=0\). Suppose \(a\) is a monomorphism in \(\clA\). As \(\clA\) is a proper abelian subcategory, there is a short triangle \(A\xrightarrow{a} C_P\xrightarrow{v} B\) in \(\clT\) with \(B\) an object in \(\clA\). As \( c_Pa=0\) and \(v\) is a weak cokernel of \(a\), there exists a morphism \(B\xrightarrow{b}{\bbS}P\) such that \( c_P=bv\). As \( c_P\) is an \(\clA\)-cover, there exists a morphism \(B\xrightarrow{\overline{b}} C_P\) such that \(b= c_P\overline{b}\). The situation is depicted in the following diagram:
    \[
        \begin{tikzcd}
            A \arrow[r, "a"] \arrow[rd, "0"'] &  C_P \arrow[r, "v"] \arrow[d, " c_P"] & B \arrow[ld, "b", dashed] \arrow[l, "\overline{b}"', dashed, bend right] \\
                 & {\bbS}P.                                           &                                                                         
        \end{tikzcd}
    \]
    We have the following equalities: \( c_P=bv= c_P\overline{b}v\). As \( c_P\) is right minimal, \(\overline{b}v\) is an automorphism. As the composition of consecutive morphisms in a triangle vanishes, we have \(\overline{b}va=0\).  It follows that \(a=0\).

    Now suppose \(a\) is arbitrary. We choose a factorisation \(a\colon A \xrightarrow{a_e} D \xrightarrow{a_m}  C_P\) in \(\clA\) where \(a_e\) is an epimorphism and \(a_m\) is a monomorphism. It suffices now to show \(a_m=0\). As \(\clA\) is a proper abelian subcategory, there is a short triangle \(A \xrightarrow{a_e} D \xrightarrow{d} \Sigma B\) with \(B\) an object in \(\clA\). As \( c_Pa_ma_e= c_Pa=0\) and \(d\) is a weak cokernel of \(a_e\), there exists a morphism \(\Sigma B \xrightarrow{b_1} {\bbS}P\) such that \( c_Pa_m=b_1d\). The situation is depicted in the following diagram:
    \[
        \begin{tikzcd}
            A \arrow[r, "a_e"] \arrow[rdd, "0"'] \arrow[rd, "a"] & D \arrow[r, "d"] \arrow[d, "a_m"] & \Sigma B \arrow[ldd, "b_1", dashed] \\
                                                     & C_P \arrow[d, "c_P"]              &                                     \\
                                                     & {\bbS}P.                                &                                    
        \end{tikzcd}
    \]
    The morphism \(b_1\) lies in the \(k\)-vector space \(\clT(\Sigma B,{\bbS}P)\). But, \(\clT(\Sigma B,{\bbS}P) \cong \D\clT(P,\Sigma B)=0\), where the isomorphism is given by Serre duality and where the equality holds as \(P\) is a projective objective in \(\clA\) (see \cite[Theorem on page~1]{dyer-exact-2005}). Hence, \(b_1=0\) and so \( c_Pa_m=0\). As \(a_m\) is a monomorphism in \(\clA\), the previous argument tells us that \(a_m=0\), as required.

    \textit{Part 2.} Dual to \textit{part 1}.
\end{proof}
\end{prop}

\begin{prop}\label{PROP: nu_P is injective and nu_I^- is projective}
    Let \(\clA\) be a category with Setup~\ref{SETUP: Base of proper abelian}. Let \(X\) be an object in \(\clA\). Then the following hold:
    \begin{itemize}
        \item[1.] If \( C_X\xrightarrow{} {\bbS}X\) is an \(\clA\)-cover in \(\clT\), then \( C_X\) is an injective object in \(\clA\).
        \item[2.] If \({\bbS}^{-1}X\xrightarrow{} E_X\) is an \(\clA\)-envelope in \(\clT\), then \(E_X\) is a projective object in \(\clA\).
    \end{itemize}
\begin{proof}
    \textit{Part 1.} Appealing to the dual of \cite[Theorem on page~1]{dyer-exact-2005}, we will show that \(\clT(A,\Sigma C_X)=0\) for each object \(A\) in \(\clA\). To this end, let \(A\) be an object in \(\clA\) and complete the morphism \( C_X\xrightarrow{} {\bbS}X\) to a short triangle \( C_X \xrightarrow{} {\bbS}X \xrightarrow{} Z\). By properties of triangulated categories, the sequence of \(k\)-vector spaces \(\clT(A,Z) \xrightarrow{} \clT(A,\Sigma C_X) \xrightarrow{} \clT(A,\Sigma {\bbS}X)\) is exact. The Triangulated Wakamatsu's Lemma \cite[Lemma~2.1]{jorgensen-auslander-2009} ensures \(\clT(A,Z)=0\). By Serre duality we have \(\clT(A,\Sigma {\bbS}X) \cong \D\clT(X,\Sigma^{-1}A)=0\), where the latter \(k\)-vector space vanishes since \(\clT(\clA,\Sigma^{-1}\clA)=0\). It follows that \(\clT(A,\Sigma C_X)=0\) by the exactness of the aforementioned sequence.
    
    \textit{Part 2.} Dual to \textit{part 1}.
\end{proof}
\end{prop}

\begin{thm}\label{THM: Nakayama functors defined and are equivalances}
    Let \(\clA\) be a category with Setup~\ref{SETUP: Base of proper abelian}. Then the following hold:
    \begin{itemize}
        \item[1.] If for each projective object \(P\) in \(\clA\) there is an \(\clA\)-cover in \(\clT\) of the form \(\nu P\xrightarrow{\alpha_P} {\bbS}P\), then the assignment \(P \mapsto \nu P\) augments to an additive functor \(\Proj\clA\xrightarrow{\nu}\Inj\clA\) such that the induced diagram
        \[
        \begin{tikzcd}
            \nu P \arrow[d, "\nu p"'] \arrow[r, "\alpha_P"] & {\bbS}P \arrow[d, "{\bbS}p"] \\
            \nu{P'} \arrow[r, "\alpha_{P'}"]                & {\bbS}P'               
        \end{tikzcd}
        \]
        is commutative for each morphism \(P \xrightarrow{p} P'\) in \(\Proj \clA\).
        
        \item[2.] If for each injective object \(I\) in \(\clA\) there is an \(\clA\)-envelope in \(\clT\) of the form \({\bbS}^{-1}I\xrightarrow{\beta_I} \nu^{-}I\), then the assignment \(I \mapsto \nu^{-} I\) augments to an additive functor \(\Inj\clA\xrightarrow{\nu^-}\Proj\clA\) such that the induced diagram
        \[
        \begin{tikzcd}
            {\bbS}^{-1}I \arrow[r, "\beta_I"] \arrow[d, "{\bbS}^{-1}i"'] & \nu^-I \arrow[d, "\nu^-{i}"] \\
            {\bbS}^{-1}I' \arrow[r, "\beta_{I'}"]                      & \nu^-{I'}                   
        \end{tikzcd}
        \]
        is commutative for each morphism \(I \xrightarrow{i} I'\) in \(\Inj \clA\).
    \end{itemize}
    Moreover, if the conditions from both \textit{part 1} and \textit{part 2} are satisfied then the functors \(\nu\) and \(\nu^-\) are mutual quasi-inverses. 
\begin{proof}
    \textit{Part 1.} The existence of a functor \(\Proj\clA\xrightarrow{}\clA\) given by the assignment \(P\mapsto \nu P\) is guaranteed by Lemma~\ref{LEM: Nakayama functor defined 0} and Proposition~\ref{PROP: strong covers and envelopes}. Proposition~\ref{PROP: nu_P is injective and nu_I^- is projective} assures that this assignment is a well-defined additive functor \(\Proj\clA \xrightarrow{\nu} \Inj\clA\).
    
    \textit{Part 2.} Dual to \textit{part 1}.

    We now show that \(\nu\) and \(\nu^-\) are mutual quasi-inverses. As both \(\nu\) and \(\nu^-\) are additive functors, it suffices to first construct a collection of isomorphisms
    \[
    \{I \xrightarrow{} \nu P \mid \text{ with } P=\nu^-I \text{ and } I \text{ an indecomposable injective object in } \clA\}
    \]
    natural in \(I\), and then to construct a collection of isomorphisms
    \[
    \{\nu^-J \xrightarrow{} Q \mid \text{ with } J=\nu Q \text{ and } Q \text{ an indecomposable projective object in } \clA\}
    \]
    natural in \(Q\). To this end, let \(I\) be an indecomposable injective object of \(\clA\) and set \(P=\nu^-I\). We have an \(\clA\)-envelope \({\bbS}^{-1}I \xrightarrow{\beta_I} P\) which is automatically a strong \(\clA\)-envelope by Proposition~\ref{PROP: strong covers and envelopes}. As \({\bbS}\) is an autoequivalence, we get a strong \({\bbS}\clA\)-envelope \(\gamma_I\colon I \cong {\bbS}{\bbS}^{-1}I \xrightarrow{{\bbS}\beta_I} {\bbS}P\). We first prove that \(\gamma_I\) is also a strong \(\clA\)-cover. 
    
    Let \(A\) be an object in \(\clA\). By Serre duality, we have the following commutative diagram:
    \begin{equation}\label{CD1}
    \begin{tikzcd}[row sep=5ex]
        {\clT({\bbS}P,{\bbS}A)} \arrow[d, "{\clT(\gamma_I,{\bbS}A)}"'] \arrow[r, "\cong"] & {\D\clT(A,{\bbS}P)} \arrow[d, "{\D\clT(A,\gamma_I)}"] \\
        {\clT(I,{\bbS}A)} \arrow[r, "\cong"]                                  & {\D\clT(A,I)}.                                  
    \end{tikzcd}
    \end{equation}
    As \(\gamma_I\) is a strong \({\bbS}\clA\)-envelope, the map \(\clT(\gamma_I,{\bbS}A)\) is an isomorphism. By the commutativity of diagram~\ref{CD1}, so too is \(\D\clT(A,\gamma_I)\) an isomorphism and therefore \(\gamma_I\) is a strong \(\clA\)-cover as required.
    
    Notice that the object \(P\) is projective in \(\clA\) using \textit{part 2} and the definition of \(P\). Furthermore, the \(\clA\)-cover \(\nu {P} \xrightarrow{\alpha_{P}} {\bbS}P\) is a strong \(\clA\)-cover by Proposition~\ref{PROP: strong covers and envelopes}. We now construct the morphism \(I \xrightarrow{} \nu P\). 
    
    By the uniqueness of (strong) \(\clA\)-covers, we get an (unique) isomorphism \(I \xrightarrow{\overline{\gamma_I}} \nu P\) such that \(\gamma_I=\alpha_P\overline{\gamma_I}\). We now show that \(\overline{\gamma_I}\) is natural. Let \(I \xrightarrow{i} J\) be a morphism between indecomposable injective objects of \(\clA\). We denote the morphism \(\nu^-I \xrightarrow{\nu^-i} \nu^-J\) in \(\clA\) by \(P \xrightarrow{p} Q\). Then we get unique morphisms \(\overline\gamma_I\) and \(\overline\gamma_J\) such that \(\gamma_I=\alpha_P\overline{\gamma_I}\) and \(\gamma_J=\alpha_Q\overline{\gamma_J}\) hold. By definition of \(\nu^-\), we have \(\beta_J {\bbS}^{-1}i=p\beta_I\) and therefore, after applying \({\bbS}\) to this equation and using the natural isomorphism \(\1 \xrightarrow{} {\bbS}{\bbS}^{-1}\), we get \(\gamma_J i={\bbS}(p)\gamma_I\). We obtain the diagram:
    \[
\begin{tikzcd}
I \arrow[rr, "\overline{\gamma_I}"] \arrow[rd, "\gamma_I"'] \arrow[dd, "i"'] &                                 & \nu P \arrow[ld, "\alpha_P"] \arrow[dd, "\nu p"] \\
                                                                             & {\bbS}P \arrow[dd, "{\bbS}p"  {yshift=16pt}] &                                                  \\
J \arrow[rr, "\overline{\gamma_J}"{xshift=-20pt}] \arrow[rd, "\gamma_J"']                  &                                 & \nu Q \arrow[ld, "\alpha_Q"]                     \\
                                                                             & {\bbS}Q.                              &                                                 
\end{tikzcd}
    \]
    As the right-facing square commutes by the definition of \(\nu P\), we have that all but the back-most square with vertices \(I, J, \nu Q\) and \(\nu P\) commute. By the application of Lemma~\ref{lem: back square commutes with strong cover}, we are done. Similarly, we construct the collection of natural isomorphisms \(\nu^-J \xrightarrow{} Q\). This proves that \(\nu\) and \(\nu^-\) are mutual quasi-inverses.
\end{proof}
\end{thm}

\begin{rem}\label{rem: can relax cover to precover}
    Because we work in a Krull-Schmidt setting, every \(\clA\)-precover can be made into an \(\clA\)-cover by removing direct summands from the source object of the \(\clA\)-precover (see \cite[Lemma~4.1]{jorgensen-auslander-2009}). As the dual situation is also true, the assumptions in Theorem~\ref{THM: Nakayama functors defined and are equivalances} may be relaxed.
\end{rem}

\begin{defn}\label{def: subcategory with a NF via aprox}
    Let \(\clA\) be a category with Setup~\ref{SETUP: Base of proper abelian} such that the following hold:
    \begin{itemize}
        \item For each projective object \(P\) in \(\clA\) there is an \(\clA\)-cover in \(\clT\) of the form \(\nu P\xrightarrow{\alpha_P} {\bbS}P\).
        \item For each injective object \(I\) in \(\clA\) there is an \(\clA\)-envelope in \(\clT\) of the form \({\bbS}^{-1}I\xrightarrow{\beta_I} \nu^{-}I\).
    \end{itemize}
    Then we call the functors \(\Proj\clA \xrightarrow{\nu} \Inj\clA\) and \(\Inj\clA \xrightarrow{\nu^-} \Proj\clA\) obtained from Theorem~\ref{THM: Nakayama functors defined and are equivalances} the \textit{Nakayama functors on \(\clA\) (obtained via approximations)} and we say that \textit{\(\clA\) has a Nakayama functor}.
\end{defn}

\begin{prop}\label{PROP: equivalent definition of nu}
    Let \(\clA\) be a category with Setup~\ref{SETUP: Base of proper abelian}. Let \(X\) and \(P\) be objects in \(\clA\) with \(P\) projective. Then the following three statements are equivalent:
    \begin{itemize}
        \item[1.] There is an \(\clA\)-cover in \(\clT\) of the form \(X \xrightarrow{\alpha_P} {\bbS}P\).
        \item[2.] There is a strong \(\clA\)-cover in \(\clT\) of the form \(X \xrightarrow{\alpha_P} {\bbS}P\).
        \item[3.] There is a natural isomorphism of the form \(\D\clA(P,-) \cong \clA(-, X)\) of functors \(\clA^{\op} \xrightarrow{} \mod_k\), i.e. the object \(X\) represents the functor \(\D\clA(P,-)\).
    \end{itemize}
    Dually, let \(Y\) and \(I\) be objects in \(\clA\) with \(I\) injective. Then the following three statements are equivalent:
    \begin{itemize}
        \item[1'.] There is an \(\clA\)-envelope in \(\clT\) of the form \({\bbS}^{-1}I \xrightarrow{} Y\).
        \item[2'.] There is a strong \(\clA\)-envelope in \(\clT\) of the form \({\bbS}^{-1}I \xrightarrow{} Y\).
        \item[3'.] There is a natural isomorphism of the form \( \D\clA(-,I) \cong \clA(Y,-)\) of functors \(\clA \xrightarrow{} \mod_k\), i.e. the object \(Y\) represents the functor \(\D\clA(-,I)\).
    \end{itemize}
\begin{proof}
    \textit{Part 1.} \((1\Leftrightarrow 2)\): This follows from Proposition~\ref{PROP: strong covers and envelopes}.
    
    \((2\Leftrightarrow 3)\): This follows from Lemma~\ref{LEM: strong covers are representables in an abelian subcat}.
    
    \textit{Part 2.} Dual to \textit{part 1}.
\end{proof}
\end{prop}

\begin{thm}\label{THM: A is a dualising vartiety}
    Let \(\clA\) be an essentially small category that has a Nakayama functor (see Definition~\ref{def: subcategory with a NF via aprox}). If \(\clA\) has enough projectives and enough injectives, then \(\clA\) is a dualising \(k\)-variety.
\begin{proof}
    We first show that the induced functor \((\Mod_{\clA})^{\op} \xrightarrow{\bbD} \Mod_{\clA^{\op}}\) restricts to a functor \((\mod_{\clA})^{\op} \xrightarrow{} \mod_{\clA^{\op}}\). To this end, let \(\clA^{\op} \xrightarrow{F} \Mod_k\) be an object in \(\mod_{\clA}\). Firstly, assume \(F=\clA(-,I)\) for \(I\) an injective object in \(\clA\). Then we have \(\bbD (F) = \bbD\clA(-,I) \cong \restr{\clT({\bbS}^{-1}I,-)}{\clA} \cong \clA(\nu^{-}I,-)\), where the first isomorphism holds since it holds pointwise by Serre duality and the second isomorphism holds as \(\beta_I\) is a strong \(\clA\)-envelope (see Proposition~\ref{PROP: strong covers and envelopes}). Hence, \(\bbD (F)\) lies in \(\mod_{\clA^{\op}}\).

    Secondly, assume \(F = \clA(-,A)\) for \(A\) an object in \(\clA\). Then choose an injective copresentation \(A \xrightarrow{} I^0 \xrightarrow{} I^1\) of \(A\). As the Yoneda embedding \(\clA \xrightarrow{} \Mod_{\clA}\) is left exact, the sequence \(\clA(-,A) \xrightarrow{} \clA(-,I^0) \xrightarrow{} \clA(-,I^1)\) is left exact. Using the exactness of \(\bbD\), the sequence \(\bbD\clA(-,I^1) \xrightarrow{} \bbD\clA(-,I^0) \xrightarrow{} \bbD\clA(-,A)=\bbD(F)\) is right exact. As \(\mod_{\clA^{\op}}\) is closed under cokernels (see \cite[Proposition on page~41]{auslander-representation-1971}) \(\bbD(F)\) lies in \(\mod_{\clA^{\op}}\).

    Lastly, assume \(F\) is arbitrary in \(\mod_{\clA}\). Then there is a right exact sequence of the form \(\clA(-,Y) \xrightarrow{} \clA(-,X) \xrightarrow{} F\). As \(\bbD\) is exact, the sequence \(\bbD(F) \xrightarrow{} \bbD\clA(-,X) \xrightarrow{} \bbD\clA(-,Y)\) is left exact. As \(\mod_{\clA^{\op}}\) is closed under kernels (see \cite[Proposition on page~41]{auslander-representation-1971}) \(\bbD(F)\) lies in \(\mod_{\clA^{\op}}\). Moreover, the fact that the functor \(\Mod_{\clA^{\op}} \xrightarrow{\bbD} (\Mod_{\clA})^{\op}\) restricts to a functor \(\mod_{\clA^{\op}} \xrightarrow{\bbD} (\mod_{\clA})^{\op}\) follows by a dual argument to the one given above.
\end{proof}
\end{thm}

\begin{rem}
    The inclusion of points 3 and 3' in Proposition~\ref{PROP: equivalent definition of nu} was inspired by ongoing work by Nan Gao, Julian K{\"u}lshammer, Sondre Kvamme and Chrysostomos Psaroudakis presented at the \textit{Homological Algebra and Representation Theory} conference (Karlovasi, Samos, 2023). They used the existence of the representations in both points 3 and 3' to define the notion of an abelian category admitting a Nakayama functor. This inspiration was instrumental in facilitating the proof of Theorem~\ref{THM: A is a dualising vartiety}.
\end{rem}

\subsection{Proper abelian length subcategories}

\begin{thm}\label{THM: enough proj iff enough inj}
    Let \(\clA\) be a category that has a Nakayama functor (see Definition~\ref{def: subcategory with a NF via aprox}) and assume each object of \(\clA\) is of finite length. Then \(\clA\) has enough projectives if and only if \(\clA\) has enough injectives.
\begin{proof}
    \((\xRightarrow[]{\text{only if}})\): Assume \(\clA\) has enough projectives and let \(A\) be a simple object in \(\clA\). By assumption, there is a nonzero epimorphism \(P \xrightarrow{p} A\) in \(\clA\) with \(P\) a projective object in \(\clA\). We first show that the composition \(\nu {P} \xrightarrow{\alpha_{P}} {\bbS}P \xrightarrow{{\bbS}p} {\bbS}A\) is nonzero. Assume the contrary and complete \(\alpha_{P}\) to a short triangle \(\nu {P} \xrightarrow{\alpha_{P}} {\bbS}P \xrightarrow{g} Z\). As \({\bbS}(p)\alpha_{P}=0\) and as \(g\) is a weak cokernel of \(\alpha_P\), we have that \({\bbS}p\) factors through a morphism in the \(k\)-vector space \(\clT(Z,{\bbS}A)\). But, \(\clT(Z,{\bbS}A) \cong \D\clT(A,Z)=0\) where the isomorphism holds by Serre duality and the latter \(k\)-vector space vanishes by the Triangulated Wakamatsu's Lemma \cite[Lemma~2.1]{jorgensen-auslander-2009}. Hence, we get \({\bbS}p=0\) and therefore \(p=0\), a contradiction.

    We now have that the composition \(\nu {P} \xrightarrow{\alpha_{P}} {\bbS}P \xrightarrow{{\bbS}p} {\bbS}A\) is nonzero. By Serre duality, the existence of a nonzero morphism \(A \xrightarrow{a} \nu {P}\) is ensured. Therefore, each simple object in \(\clA\) admits a nonzero morphism to an injective object in \(\clA\) (see Proposition~\ref{PROP: nu_P is injective and nu_I^- is projective}) and by Lemma~\ref{LEM: finite length categories, simples and enough injectives/projectives}, \(\clA\) has enough injectives.

    \((\xLeftarrow[]{\text{if}})\): Dual to the previous argument.
\end{proof}
\end{thm}

\section{Auslander-Reiten translates on proper abelian subcategories}
\subsection{The construction and fundamental properties}\label{SUBSEC: the construction and fundamental properties}

\begin{setup}\label{SETUP:Proper abelian with all info.}
    Throughout this section, fix a category \(\clA\) that has a Nakayama functor (see Definition~\ref{def: subcategory with a NF via aprox}) and assume that each object in \(\clA\) has a projective cover and an injective envelope.
\end{setup}

    For each indecomposable non-projective object \(A\) in \(\clA\) and for each indecomposable non-injective object \(B\) in \(\clA\) we fix
    \[
    P_1 \xrightarrow{p_1} P_0 \xrightarrow{p_0} A \text{\qquad and \qquad} B \xrightarrow{i^0} I^0 \xrightarrow{i^1} I^1
    \]
    to be a minimal projective presentation of \(A\) and a minimal injective copresentation of \(B\), respectively. For each indecomposable projective object \(P\) in \(\clA\) and for each indecomposable injective object \(I\) in \(\clA\), we fix
    \[
    0 \xrightarrow{} P \xrightarrow{1} P \text{\qquad and \qquad} I \xrightarrow{1} I \xrightarrow{} 0
    \]
    to be a minimal projective presentation of \(P\) and a minimal injective copresentation of \(I\), respectively. Furthermore, we extend to general objects of \(\clA\) by taking direct sums.

\begin{defn}\label{DEF: Auslander-Reiten translates}
     Let \(C\) be an object in \(\clA\). We define the \textit{Auslander-Reiten translate \(\tau C\) of \(C\)} as \(\tau C = \Ker(\nu Q_1 \xrightarrow{\nu q_1} \nu Q_0)\), where \(Q_1 \xrightarrow{q_1} Q_0 \xrightarrow{q_0} C\) is the fixed minimal projective presentation of \(C\). Similarly, we define the \textit{(inverse) Auslander-Reiten translate \(\tau^- C\) of \(C\)} as \(\tau^- C = \Coker(\nu^- J^0 \xrightarrow{\nu^-{j^1}} \nu^- J^1)\), where \(C \xrightarrow{j^0} J^0 \xrightarrow{j^1} J^1\) is the minimal injective copresentation of \(C\).
\end{defn}
    
\begin{prop}\label{PROP: Properties of tau.}
    Let \(\clA\) be a category with Setup \ref{SETUP:Proper abelian with all info.} and let \(A\) and \(A'\) be indecomposable objects in \(\clA\). Then the following hold:
    \begin{itemize}
        \item[1.] \(\tau(A\oplus A') \cong \tau A \oplus\tau A'\).
        \item[2.] The object \(A\) is projective in \(\clA\) \(\iff\) \(\tau A=0\).
        \item[3.] If \(A\) is a non-projective object in \(\clA\), then \(\tau A\) is a non-injective object in \(\clA\).
        \item[4.] If \(A\) is a non-projective object in \(\clA\), then the sequence \(\tau A \xrightarrow{k_A} \nu P_1 \xrightarrow{\nu p_1} \nu P_0\) is a minimal injective copresentation of \(\tau A\).
        \item[5.] If \(A\) is a non-projective object in \(\clA\), then \(\tau^-\tau A\cong A\).
        \item[6.] If \(A\) is a non-projective object in \(\clA\), then \(\tau A\) is indecomposable.
        \item[7.] If \(A\) and \(A'\) are non-projective objects in \(\clA\), then \(A\cong A' \iff \tau A \cong \tau A'\).
    \end{itemize}
    \begin{proof}
    \textit{Part 1.} We have minimal projective presentations \(P_1 \xrightarrow{p_1} P_0 \xrightarrow{p_0} A\) and \(P'_1 \xrightarrow{p'_1} P'_0 \xrightarrow{p'_0} A'\) of \(A\) and \(A'\) respectively and therefore, \(P_1\oplus P_1' \xrightarrow{p_1\oplus p'_1} P_0\oplus P_0' \xrightarrow{p_0\oplus p_0'} A\oplus A'\) is the fixed minimal projective presentation of \(A\oplus A'\). By the application of \(\nu\), we arrive at the following solid commutative diagram:
    \[
    \begin{tikzcd}[column sep=10ex]
        0 \arrow[r] & \tau(A\oplus A') \arrow[r] \arrow[d, "\cong", dashed] & \nu(P_1\oplus P_1') \arrow[r, "\nu(p_1\oplus p_1')"] \arrow[d, "\cong"] & \nu(P_0\oplus P_0') \arrow[d, "\cong"] \\
        0 \arrow[r] & \tau A\oplus \tau A' \arrow[r]                        & \nu P_1 \oplus \nu P_1' \arrow[r, "\nu p_1 \oplus \nu p_1'"]            & \nu P_0 \oplus \nu P_0'
    \end{tikzcd}
    \]
    where the vertical solid morphisms are the canonical isomorphisms induced by the additivity of \(\nu\) (see Theorem \ref{THM: Nakayama functors defined and are equivalances}). The dashed morphism exists by the universal property of the kernel of \(\nu p_1 \oplus \nu p_1'\) and the Five Lemma ensures it is an isomorphism.
    
    \textit{Part 2.} \((\Rightarrow)\): Assume \(A\) is a projective object in \(\clA\). Then \(0 \xrightarrow{} A \xrightarrow{1} A\) is the fixed minimal projective presentation of \(A\). As \(0=\nu 0 \xrightarrow{} \nu A\) is a monomorphism, it follows that \(\tau A=0\).
    
    \((\Leftarrow)\): Assume \(\tau A=0\). Then \(\nu {P_1} \xrightarrow{\nu {p_1}} \nu {P_0}\) is a monomorphism in \(\clA\) with \(\nu {P_1}\) an injective object in \(\clA\). Therefore, \(\nu p_1\) is a section. As \(\nu\) is an equivalence, the morphism \(p_1\) must also be a section. Hence, we have \(P_0 \cong P_1 \oplus A\) and therefore, \(A\) is projective.
    
    \textit{Part 3.} Assume that \(A\) is a non-projective object in \(\clA\). By Lemma~\ref{LEM: A specific precover and preenvelope}, there is a nonzero \(\clA\)-precover of the form \(\tau A \xrightarrow{\theta} \Sigma^{-1}\bbS A\) (notice that \(\tau A\) is \(X_A\) in the statement of Lemma~\ref{LEM: A specific precover and preenvelope}). Showing that \(\tau A\) is not an injective object in \(\clA\) amounts to finding an object \(C\) in \(\clA\) with \(\clT(C,\Sigma\tau A)\neq 0\) (see dual of \cite[Theorem on page~1]{dyer-exact-2005}). Choosing \(C=A\) and using Serre duality, we have \(\clT(A,\Sigma\tau A) \cong \D\clT(\tau A,\Sigma^{-1}{\bbS}A)\), where the last \(k\)-vector space is non-vanishing since the nonzero morphism \(\theta\) lies in \(\clT(\tau A,\Sigma^{-1}{\bbS}A)\).
    
    \textit{Part 4.} It is clear the sequence is an injective copresentation of \(\tau A\). As \(\nu\) is an additive equivalence (see Theorem \ref{THM: Nakayama functors defined and are equivalances}), Lemma \ref{LEM: isomorphism between radical maps induced by an additive equivalence} tells us that \(\nu p_1\) is a radical morphism as, by \cite[Proposition~3.10]{krause-krull-2015}, \(p_1\) is a radical morphism. The fact that \(k_A\) is an injective envelope follows from the dual of \cite[Proposition~3.10]{krause-krull-2015}. As \(\nu(p_1) k_A=0\), the universal property of the cokernel of \(k_A\) ensures the existence of a unique morphism \(\Omega^{-1}\tau A \xrightarrow{\iota_1} \nu P_0\) making the following diagram commutative:
\[
    \begin{tikzcd}
    0 \arrow[r] & \tau A \arrow[r, "k_A"] & \nu P_1 \arrow[r, "\nu p_1"] \arrow[d, "\pi_1"] & \nu P_0 \\
            &                         & \Omega^{-1}\tau A. \arrow[ru, "\iota_1"', dashed]                 &                        &  
    \end{tikzcd}
\]
    Here, \(\Omega^{-1}\tau A\) is the cosyzygy of \(\tau A\) and \(\pi_1\) is the cokernel of \(k_A\). By diagram chasing, \(\iota_1\) is a monomorphism, so it suffices to show that \(\iota_1\) is an injective envelope.
    
    To this end, we choose an injective envelope \(I\) of \(\Coker \nu p_1\) extending the injective copresentation of \(\tau A\):
    \[
    \begin{tikzcd}
    0 \arrow[r] & \tau A \arrow[r, "k_A"] & \nu P_1 \arrow[r, "\nu p_1"] \arrow[d, "\pi_1"] & \nu P_0 \arrow[r, "f"] & I \\
            &                         & \Omega^{-1}\tau A. \arrow[ru, "\iota_1"']        &                        &  
    \end{tikzcd}
    \]
    By the dual of \cite[Proposition~3.10]{krause-krull-2015}, Lemma~\ref{LEM: isomorphism between radical maps induced by an additive equivalence} and the fact that \(\nu\) and \(\nu^-\) are mutual quasi-inverses (see Theorem~\ref{THM: Nakayama functors defined and are equivalances}), it suffices to show that the morphism \(g\colon P_0 \xrightarrow{\cong} \nu^-\nu P_0 \xrightarrow{\nu^-f} \nu^- I\) is a radical morphism. As \(gp_1=0\) and \(p_0\) is the cokernel of \(p_1\), there exists a unique morphism \(A \xrightarrow{a} \nu^- I\) such that the following diagram is commutative:
    \[
    \begin{tikzcd}
    P_1 \arrow[r, "p_1"] & P_0 \arrow[r] \arrow[r, "p_0"] \arrow[rd, "g"'] &    A \arrow[r] \arrow[d, "a"] & 0 \\
                     &                                                 & \nu^-I.                    &  
\end{tikzcd}
    \]
    But, \(p_0\) must be a radical morphism as \(A\) is a non-projective indecomposable object. Hence, \(g=ap_0\) is a radical morphism.
    
    \textit{Part 5.} Assume \(A\) is a non-projective object in \(\clA\). By \textit{part 4}, \(\tau A \xrightarrow{} \nu P_1 \xrightarrow{\nu p_1} \nu P_0\) is a minimal injective copresentation of \(\tau A\). Therefore, as minimal injective copresentations are unique up to isomorphism, we get the following commutative diagram:
    \[
    \begin{tikzcd}
        0 \arrow[r] & \tau A \arrow[r] & \nu P_1 \arrow[r, "\nu p_1"] \arrow[d, "\cong"] & \nu P_0 \arrow[d, "\cong"] \\
        0 \arrow[r] & \tau A \arrow[r] & I^0 \arrow[r, "i^0"]                            & I^1                       
    \end{tikzcd}
    \]
    where the bottom row is our fixed injective copresentation of \(\tau A\). Applying \(\nu^-\), we get the following solid commutative diagram:
    \[
    \begin{tikzcd}
        P_1 \arrow[r, "p_1"] \arrow[d, "\cong"]                   & P_0 \arrow[r, "p_0"] \arrow[d, "\cong"]                   & A \arrow[r] \arrow[d, "\cong", dashed] & 0 \\
        \nu^-\nu P_1 \arrow[r, "\nu^-\nu p_1"] \arrow[d, "\cong"] & \nu^-\nu P_0 \arrow[r, "\nu^-\nu p_0"] \arrow[d, "\cong"] & C \arrow[r] \arrow[d, "\cong", dashed] & 0 \\
        \nu^- I^0 \arrow[r, "\nu^-i^1"]                           & \nu^- I^1 \arrow[r]                                       & \tau^-\tau A \arrow[r]                 & 0
    \end{tikzcd}
    \]
    where \(C\) is the cokernel of \(\nu^-\nu p_1\). Then the dashed morphisms exist by the universal property of the cokernel and the Five Lemma ensures that they are isomorphisms.
    
    \textit{Part 6.} Suppose \(A\) is a non-projective object in \(\clA\) and assume, for a contradiction, we have \(\tau A\cong X\oplus Y\), for nonzero objects \(X\) and \(Y\) in \(\clA\). We first show that both \(X\) and \(Y\) are non-injective objects in \(\clA\).

    For a contradiction, assume \(Y\) is a nonzero injective object in \(\clA\). As \(s\colon Y \xrightarrow{} \tau A \xrightarrow{k_A} \nu P_1\) is a monomorphism, where \(Y \xrightarrow{} \tau A\) is the canonical inclusion morphism, \(s\) must be a section. As \(\nu(p_1) s=0\), the object \(Y\) is isomorphic to a direct summand of \(\nu P_1\) that lies in the kernel of \(\nu p_1\). By applying \(\nu^-\), we get a section \(\nu^-s\) with \(\nu^-(\nu (p_1)s)=0\). Similarly, the object \(\nu^-Y\) is isomorphic to a direct summand of \(\nu^-\nu P_1\) that lies in the kernel of \(\nu^-\nu p_1\). As \(\nu\) and \(\nu^-\) are mutual quasi-inverses (see Theorem~\ref{THM: Nakayama functors defined and are equivalances}), there is a nonzero direct summand of \(P_1\) contained in the kernel of \(p_1\). Using that the projective cover \(P_1 \xrightarrow{} \Omega A\) is unique up to isomorphism, we get a contradiction to the dual of \cite[Corollary~1.4]{krause-minimal-1998}.
    
    Now, we have \(A\cong\tau^-\tau A \cong \tau^-X \oplus \tau^-Y\) by \textit{part 5} and the dual of \textit{part 1}. By the dual of \textit{part 2}, both objects \(\tau^-X\) and \(\tau^-Y\) are nonzero, a contradiction as \(A\) is indecomposable.

    \textit{Part 7.} \((\Rightarrow)\): Assuming \(A\cong A'\), the result follows by the uniqueness of minimal projective presentations, the universal property of kernels and the Five Lemma.

    \((\Leftarrow)\): Conversely, assume \(\tau A\cong \tau A'\). By \textit{part 5}, we have \(A\cong\tau^-\tau A\) and \(A'\cong\tau^-\tau A'\). Then by \textit{part 6} and \textit{part 3}, \(\tau A\) and \(\tau A'\) are indecomposable non-injective objects in \(\clA\). Therefore, the dual argument of the ``\(\Rightarrow\)'' implication of \textit{part 7} gives us \(\tau^-\tau A\cong \tau^-\tau A'\).
    \end{proof}
\end{prop}

    We state the dual version of the above result for completeness.

\begin{prop}\label{PROP: properties of inverse tau}
    Let \(\clA\) be a category with Setup \ref{SETUP:Proper abelian with all info.} and let \(B\) and \(B'\) be indecomposable objects in \(\clA\). Then the following hold:
    \begin{itemize}
        \item[1.] \(\tau^-(B\oplus B') \cong \tau^-B\oplus \tau^-B'\).
        \item[2.] The object \(B\) is injective in \(\clA\) \(\iff\) \(\tau^- B=0\).
        \item[3.] If \(B\) is a non-injective object in \(\clA\), then \(\tau^- B\) is a non-projective object in \(\clA\).
        \item[4.] If \(B\) is a non-injective object in \(\clA\), then the sequence \(\nu^- I^0 \xrightarrow{} \nu^- I^1 \xrightarrow{} \tau^- B\) is a minimal projective presentation of \(\tau^- B\).
        \item[5.] If \(B\) is a non-injective object in \(\clA\), then \(\tau\tau^- B\cong B\).
        \item[6.] If \(B\) is a non-injective object in \(\clA\), then \(\tau^- B\) is indecomposable.
        \item[7.] If \(B\) and \(B'\) are non-injective objects in \(\clA\), then \(B\cong B'\) \(\iff\) \(\tau^- B \cong \tau^- B'\).
    \end{itemize}
\end{prop}

\begin{cor}
    Let \(\clA\) be a category with Setup \ref{SETUP:Proper abelian with all info.}. There is a bijective correspondence between the following:
    \begin{itemize}
        \item[1.] Isomorphism classes of indecomposable non-projective objects \(A\) in \(\clA\).
        \item[2.] Isomorphism classes of indecomposable non-injective objects \(B\) in \(\clA\).
    \end{itemize}
    The bijective correspondence is given by the assignments \(A \mapsto \tau A\) and \(B \mapsto \tau^{-}B\).
\end{cor}

\subsection{The existence of Auslander-Reiten sequences in proper abelian subcategories}

\begin{rem}\label{REM: removing essentially small triangulated category}
    In \cite[Theorem~3.1]{jorgensen-auslander-2009}, the triangulated categories are assumed to be essentially small. This assumption can be dropped due to \cite[Proposition~5.15]{iyama-auslander-2024}.
\end{rem}

\begin{thm}\label{THM: Existence of AR sequences in subcateories}
    Let \(\clA\) be a category with Setup \ref{SETUP:Proper abelian with all info.}. Then the following hold:
    \begin{itemize}
        \item [1.] For each indecomposable non-projective object \(C\) in \(\clA\), there exists an Auslander-Reiten sequence in \(\clA\) of the form
            \[
            0 \xrightarrow{} \tau C \xrightarrow{} E_C \xrightarrow{} C \xrightarrow{} 0.
            \]
        \item [2.] For each indecomposable non-injective object \(A\) in \(\clA\), there exists an Auslander-Reiten sequence in \(\clA\) of the form
            \[
            0 \xrightarrow{} A \xrightarrow{} F_A \xrightarrow{} \tau^-A \xrightarrow{} 0.
            \]
    \end{itemize}
\begin{proof}
    \textit{Part 1.} By Lemma \ref{LEM: A specific precover and preenvelope}, there exists a nonzero \(\clA\)-precover of the form \(\tau C \xrightarrow{} \Sigma^{-1}{\bbS}C\). But, \(\tau C\) is an indecomposable object in \(\clA\) by Proposition \ref{PROP: Properties of tau.} and hence by the dual of \cite[Lemma~2.4]{krause-auslander-2000}, \(\tau C \xrightarrow{} \Sigma^{-1}{\bbS}C\) is an \(\clA\)-cover. The existence of the required Auslander-Reiten sequence is ensured by \cite[Theorem~3.1]{jorgensen-auslander-2009} (see also Remark \ref{REM: removing essentially small triangulated category}).

    \textit{Part 2.} Dual to \textit{part 1}.
\end{proof}
\end{thm}

\section{Application to the module category of a finite dimensional algebra}

    In this section, we recover the standard Nakayama functors associated with the category of finite dimensional modules over a finite dimensional algebra and give a new proof of the existence of Auslander-Reiten sequences in such a category. The following setup applies throughout this section.

\begin{setup}
    Let \(A\) be a finite dimensional \(k\)-algebra and let \(\mod_A\) denote the category of finite dimensional right modules over \(A\). We also let \(\proj_A\) denote the full subcategory of finitely generated projective \(A\)-modules and let \(\inj_A\) denote the full subcategory of finitely generated injective \(A\)-modules.
\end{setup}

\begin{rem}\label{REM: removing algebraically closed field}
    In \cite[Theorem~I.2.4]{reiten-noetherian-2002}, the ground field is assumed to be algebraically closed. This assumption can be dropped due to \cite[Theorem~3.6]{iyama-auslander-2024}.
\end{rem}

    By \cite[Proposition~I.1.4]{reiten-noetherian-2002} and \cite[Proposition~I.2.3]{reiten-noetherian-2002} (see also Remark~\ref{REM: removing algebraically closed field}), we can associate to the Serre functor \(\bbS\) and the collection of natural isomorphisms 
    \[
    \{\clT(X,Y) \xrightarrow{\eta_{X,Y}} \D\clT(Y,\bbS X) \mid\text{for } X \text{ and }Y \text{ objects in }\clT\},
    \]
    a collection of linear forms
    \[
    \left\{\clT(Z,\bbS Z) \xrightarrow{\eta_Z} k \bigm\vert \text{ for } Z \text{ an object in }\clT \text{ and } \eta_Z\coloneq\eta_{Z,Z}(1_Z)\right\},
    \]
    satisfying \(\eta_Z(h_Z)\neq 0\) for each connecting homomorphism \(h_Z\) in an Auslander-Reiten triangle of the form \(
    \Sigma^{-1}{\bbS}Z \xrightarrow{} Y \xrightarrow{} Z \xrightarrow{h_Z} {\bbS}Z\). Then for each morphism \(X\xrightarrow{f} Y\) we have \(\eta_{X,Y}(f)(-)=\eta_X(-\circ f)\) as linear forms on \(\clT(Y,\bbS X)\). It will be also useful to note that, for \(X\) and \(Y\) objects in \(\clT\), the morphism \(\clT(X, \bbS Y) \xrightarrow{\operatorname{ev}} \D^2\clT(X,\bbS Y) \xrightarrow{\D(\eta_{Y,X})} \D\clT(Y,X)\) is given by \(g \mapsto \eta_Y(g \circ -)\), where \(\operatorname{ev}\) is the  evaluation map isomorphism.

\begin{rem}\label{Rem: Nakayama duality}
    Similarly to Serre duality, the standard Nakayama functor \(\operatorname{N}(-)= - \otimes_A\D(A)\) on \(\mod_A\) exhibits a duality as follows: For each finitely generated \(A\)-module \(M\) and for each finitely generated projective \(A\)-module \(P\), we have the following composition of natural isomorphisms:
    \[
    \Hom_A(P,M) \xrightarrow{} \D(P\otimes_A \D(M)) \xrightarrow{} \D(P\otimes_A \Hom_A(M,\D(A))) \xrightarrow{}  \D(\Hom_A(M, \N(P)).
    \]
    Here, the first isomorphism uses the evaluation isomorphism \(M\cong\D^2(M)\) and Tensor-Hom Adjunction, the second isomorphism uses the canonical isomorphism \(M \cong M \otimes_A A\) and the Tensor-Hom Adjunction and the third isomorphism follows by the tensor evaluation isomorphism (a reference for these isomorphisms can be found in \cite[Section~I.1]{christensen-derived-2024}). We will denote this composition by
    \[
    \Hom_A(P,M) \xrightarrow{\varepsilon_{P,M}} \D\Hom_A(M,\N(P)).
    \]
    The isomorphism \(\varepsilon_{P.M}\) is natural in both \(P\) and \(M\). We will refer to this as \textit{Nakayama duality}.
\end{rem}

\subsection{Recovering the standard Nakayama functors on the module category}\label{section: recovering the usual Nakayama functor on modA}
    
\begin{thm}\label{THM: Nakayama functors in modA}
    Let \(\mod_A\) be a \(k\)-linear subcategory of \(\clT\). Let \(\operatorname{N}(-)= - \otimes_A\D(A)\) and let \(\operatorname{N}^-(-)=\Hom_A(\D(A),-)\) be the standard Nakayama functors on \(\mod_A\). Then the following hold:
    \begin{itemize}
        \item[1.] For each indecomposable finitely generated projective \(A\)-module \(P\), there is a strong \(\mod_A\)-cover of the form \(\N(P) \xrightarrow{\alpha_P} {\bbS}P\) such that the diagram
        \[
        \begin{tikzcd}
            \N(P) \arrow[r, "\alpha_P"] \arrow[d, "\N(p)"'] & \bbS P \arrow[d, "\bbS p"] \\
            \N(Q) \arrow[r, "\alpha_Q"]                     & \bbS Q                
        \end{tikzcd}
        \]
        is commutative for each \(A\)-module homomorphism \(P \xrightarrow{p} Q\) between indecomposable finitely generated projective \(A\)-modules. Moreover, the system of strong \(\mod_A\)-covers \(\{\N(P) \xrightarrow{\alpha_P} \bbS P\}\) given above gives rise to a functor \(\proj_A \xrightarrow{\nu} \inj_A\) by Theorem~\ref{THM: Nakayama functors defined and are equivalances}. The restricted Nakayama functor \(\proj \xrightarrow{\operatorname{N}} \inj_A\) is naturally isomorphic to \(\nu\).
        
        \item[2.] For each indecomposable finitely generated injective \(A\)-module \(I\), there is a strong \(\mod_A\)-envelope of the form \({\bbS}^{-1}I \xrightarrow{} \N^{-}(I)\) such that the diagram
        \[
        \begin{tikzcd}
            \S^{-1}I \arrow[r] \arrow[d, "\S^{-1}i"'] & \N^{-}(I) \arrow[d, "\N^{-}(i)"] \\
            \S^{-1}J \arrow[r]                        & \N^{-}(J)                         
        \end{tikzcd}
        \]
        is commutative for each \(A\)-module homomorphism \(I \xrightarrow{i} J\) between indecomposable finitely generated injective \(A\)-modules. Moreover, the system of strong \(\mod_A\)-envelopes \(\{\S^{-1}I \xrightarrow{} \N^{-1}(I)\}\) given above give rise to a functor \(\inj_A \xrightarrow{\nu^-} \proj_A\) by Theorem~\ref{THM: Nakayama functors defined and are equivalances}. The restricted functor \(\inj_A \xrightarrow{\operatorname{N}^-} \proj_A\) is naturally isomorphic to \(\nu^-\).
    \end{itemize}
\begin{proof}
    \textit{Part 1.}
    \paragraph{The construction of \(\alpha_P\).} 
        Using Remark~\ref{Rem: Nakayama duality}, we have a natural isomorphism:
        \[
        \gamma_P^{(-)} \colon \Hom_A(-,\N(P)) \xrightarrow{\D(\varepsilon_{P,-})\circ \operatorname{ev}} \D\Hom_A(P,-) \xrightarrow{(\D(\eta_{P,-}) \circ \operatorname{ev})^{-1}} \clT(-,\bbS P).
        \]
        Let \(\alpha_P\coloneq\gamma_P^{\N(P)}(1_{\N(P)})\), where \(\gamma_P^{\N(P)}\) is the component of the natural isomorphism \(\gamma_P\) at \(\N(P)\).

    \paragraph{\(\alpha_P\) is a strong \(\mod_A\)-cover.}
        By definition, the morphism \(\alpha_P\) is a strong \(\mod_A\)-cover if and only if for each finitely generated \(A\)-module \(M\), the \(k\)-linear map
        \[
        \Hom_A(M,\N(P)) \xrightarrow{\clT(M,\alpha_P)} \clT(M,\bbS P)
        \]
        is an isomorphism. As \(\gamma_P^{(-)}\) is a natural isomorphism, it, therefore, suffices to show for each finitely generated \(A\)-module \(M\), the invertible component \(\gamma_P^M\) coincides with \(\Hom_A(M,\N(P)) \xrightarrow{\clT(M,\alpha_P)} \clT(M,\bbS P)\). To this end, let \(M\xrightarrow{m} \N(P)\) be an \(A\)-module homomorphism. By naturality, the following diagram is commutative:
        \[
        \begin{tikzcd}
            {\Hom_A(\N(P),\N(P))} \arrow[r, "\gamma_P^{\N(P)}"] \arrow[d, "{\Hom_A(m,\N(P))}"'] & {\clT(\N(P),\bbS P)} \arrow[d, "{\clT(m,\bbS P)}"] \\
            {\Hom_A(M,\N(P))} \arrow[r, "\gamma_P^M"]                                              & {\clT(M,\bbS P)}.                              
        \end{tikzcd}
        \]
        Chasing the identity \(1_{\N(P)}\) through this diagram yields the equality
        \[
        \gamma_P^M(m)=\clT(M,\alpha_P)(m).
        \]
        
    \paragraph{The commutativity of the induced diagram.}
        Let \(P \xrightarrow{p} Q\) be a morphism between indecomposable finitely generated projective \(A\)-modules. We want to show that the diagram
        \begin{equation}\label{Diagram:commutative square serre nak}
        \begin{tikzcd}
            \N(P) \arrow[r, "\alpha_P"] \arrow[d, "\N(p)"'] & \bbS P \arrow[d, "\bbS p"] \\
            \N(Q) \arrow[r, "\alpha_Q"]                     & \bbS Q                
        \end{tikzcd}       
        \end{equation}
        is commutative. Using the naturality of \(\gamma\) in both components, the following diagram is commutative:
        \begin{equation}\label{diagram: big comm of sere nak}
        \begin{tikzcd}
            {\Hom_A(\N(P),\N(P))} \arrow[d, "{\Hom_A(\N(P),\N(p))}"'] \arrow[r, "\cong"] & {\D\Hom_A(P,\N(P))} \arrow[d] \arrow[r, "\cong"] & {\clT(\N(P),\bbS P)} \arrow[d, "{\clT(\N(P),\bbS p)}"]  \\
            {\Hom_A(\N(P),\N(Q))} \arrow[r, "\cong"]                                       & {\D\Hom_A(Q,\N(P))} \arrow[r, "\cong"]           & {\clT(\N(P),\bbS Q)}                                  \\
            {\Hom_A(\N(Q),\N(Q))} \arrow[u, "{\Hom_A(\N(p),\N(Q))}"] \arrow[r, "\cong"]  & {\D\Hom_A(Q,\N(Q))} \arrow[u] \arrow[r, "\cong"] & {\clT(\N(Q),\bbS Q)}. \arrow[u, "{\clT(\N(p),\bbS Q)}"']
        \end{tikzcd}
        \end{equation}
        Commutativity of (\ref{Diagram:commutative square serre nak}) is equivalent to the validity of the equation \[\clT(\N(P),\bbS p)(\alpha_P)=\clT(\N(p),\bbS Q)(\alpha_Q)\] and this, by commutativity of (\ref{diagram: big comm of sere nak}), is equivalent to the validity of the equation \[\Hom_A(\N(P),\N(p))(1_{\N(P)})=\Hom_A(\N(p),\N(Q))(1_{\N(Q)}),\] which indeed holds.

    \paragraph{The construction of \(\nu\) is naturally isomorphic to \(\operatorname{N}\).} As any finitely generated projective \(A\)-module \(P\) is a finite direct sum of indecomposable finitely generated projective \(A\)-modules and as the finite direct sum of strong covers is again a strong cover, we have strong \(\mod_A\)-covers of the form \(\N(P) \xrightarrow{} \bbS P\), where \(P\) is a finitely generated projective \(A\)-module. Theorem~\ref{THM: Nakayama functors defined and are equivalances} therefore. gives rise to a functor \(\proj_A \xrightarrow{\nu} \inj_A\). The fact that the functors \(\nu\) and the restricted Nakayama functor \(\proj_A \xrightarrow{\operatorname{N}} \inj_A\) are naturally isomorphic follows from Lemma~\ref{REM: any choice of strong covers yield nat isom fucntors}.
        
    \textit{Part 2.} This follows by a similar dual argument to \textit{part 1}.
\end{proof}
\end{thm}

\subsection{Recovering the existence of Auslander-Reiten sequences in module categories}

\begin{prop}\label{Prop: certain covers and envelops in modA}
    Let \(\mod_A\) be an extension closed proper abelian \(k\)-linear subcategory of \(\clT\) and assume \(\clT(\mod_A,\Sigma^{-1}\mod_A)=0\). Let \(\operatorname{N}(-)= - \otimes_A\D(A)\) and let \(\operatorname{N}^-(-)=\Hom_A(\D(A),-)\) be the standard Nakayama functors on \(\mod_A\) and let \(\operatorname{t}\) and \(\operatorname{t}^-\) denote the standard Auslander-Reiten translates in \(\mod_A\). Then the following hold:
    \begin{itemize}
        \item[1.] If \(P_1 \xrightarrow{p_1} P_0 \xrightarrow{p_0} M\) is a projective presentation of a non-projective indecomposable finitely generated \(A\)-module \(M\), then there is a \(\mod_A\)-cover of the form \(\operatorname{t}(M) \xrightarrow{} \Sigma^{-1}{\bbS}M\) in \(\clT\), where \(\operatorname{t}(M) \xrightarrow{} \N(P_1)\) is the kernel of the induced \(A\)-module homomorphism \(\N(P_1) \xrightarrow{} \N(P_0)\) in \(\mod_A\).
        
        \item[2.] If \(L \xrightarrow{} I^0 \xrightarrow{} I^1\) is an injective copresentation of a non-injective indecomposable finitely generated \(A\)-module \(L\), then there is a \(\mod_A\)-envelope of the form \(\Sigma {\bbS}^{-1}L \xrightarrow{}\operatorname{t}^-(L)\) in \(\clT\), where \(\N^-(I^1) \xrightarrow{} \operatorname{t}^-(L) \) is the cokernel of the induced \(A\)-module homomorphism \(\N(I^0) \xrightarrow{} \N(I^1)\) in \(\mod_A\).
    \end{itemize}
\begin{proof}
    \textit{Part 1.} By Theorem~\ref{THM: Nakayama functors in modA} and Lemma~\ref{LEM: A specific precover and preenvelope}, there is a nonzero \(\mod_A\)-precover of the form \(\tau M \xrightarrow{\theta} \Sigma^{-1}{\bbS}M\), where \(\tau M \xrightarrow{} \nu P_1\) is the kernel of the induced morphism \(\nu P_1 \xrightarrow{\nu p_1} \nu P_0\) in \(\mod_A\). We have left exact sequences \(\tau M \xrightarrow{} \nu P_1 \xrightarrow{\nu p_1} \nu P_0\) and \(\operatorname{t}(M) \xrightarrow{} \N(P_1) \xrightarrow{\N(p_1)} \N(P_0)\). By the equivalence \(\nu\cong\N\) in Theorem~\ref{THM: Nakayama functors in modA} and by the universal property of the kernel of \(\N(p_1)\), we get a commutative diagram
    \[
    \begin{tikzcd}
        0 \arrow[r] & \tau M \arrow[r] \arrow[d, dashed, "f"] & \nu P_1 \arrow[r, "\nu p_1"] \arrow[d, "\cong", dashed] & \nu P_1 \arrow[d, "\cong", dashed] \\
        0 \arrow[r] & \operatorname{t}(M) \arrow[r]   & \N(P_1) \arrow[r, "\N(p_1)"]       & \N(P_0),                          
    \end{tikzcd}
    \]
    where the right two vertical \(A\)-module homomorphisms are isomorphisms. The Five Lemma ensures that the leftmost \(A\)-module homomorphism is also an isomorphism. Hence, \(\operatorname{t}(M) \xrightarrow{f^{-1}} \tau M \xrightarrow{\theta} \Sigma^{-1}{\bbS}M\) is a \(\mod_A\)-precover. Finally, the minimality of the nonzero \(\mod_A\)-precover \(\theta f^{-1}\) follows as \(\operatorname{t}(M)\) is indecomposable (see the dual of \cite[Lemma~2.4]{krause-auslander-2000}).

    \textit{Part 2.} Dual to \textit{part 1}.
\end{proof}
\end{prop}

\begin{prop}\label{THM: AR sequences in modA}
    Let \(\mod_A\) be an extension closed proper abelian \(k\)-linear subcategory of \(\clT\) and assume \(\clT(\mod_A,\Sigma^{-1}\mod_A)=0\). Let \(\operatorname{t}\) and \(\operatorname{t}^-\) denote the standard Auslander-Reiten translates in \(\mod_A\). Then the following hold:
    \begin{itemize}
        \item [1.] For each indecomposable finitely generated non-projective \(A\)-module \(M\), there is an  \linebreak Auslander-Reiten sequence in \(\mod_A\) of the form
            \[
            0 \xrightarrow{} \operatorname{t}(M) \xrightarrow{} E_M \xrightarrow{}M \xrightarrow{} 0.
            \]
        \item [2.] For each indecomposable finitely generated non-injective \(A\)-module \(L\), there is an Auslander-Reiten sequence in \(\mod_A\) of the form
            \[
            0 \xrightarrow{} L \xrightarrow{} F_L \xrightarrow{} \operatorname{t}^-(L) \xrightarrow{} 0.
            \]
    \end{itemize}
\begin{proof}
    \textit{Part 1.} By Proposition~\ref{Prop: certain covers and envelops in modA}, there exists a \(\mod_A\)-cover of the form \(\operatorname{t}(M) \xrightarrow{}\Sigma^{-1}{\bbS}M\). The existence of the required Auslander-Reiten sequence is ensured by \cite[Theorem~3.1]{jorgensen-auslander-2009} (see Remark~\ref{REM: removing essentially small triangulated category}).

    \textit{Part 2.} Dual to \textit{part 1}.
\end{proof}
\end{prop}

    The following were defined in \cite[page~351-352]{hughes-trivial-1983} and \cite[Chapter~II.2.1]{happel-triangulated-1988} (see also \cite{happel-derived-1987}): let \(\widehat{A}\) denote the repetitive algebra of \(A\), let \(\mod_{\widehat{A}}\) denote the module category of finitely generated modules over \(\widehat{A}\) and let \(\underline{\mod}_{\widehat{A}}\) denote the stable module category of \(\widehat{A}\). We now recover the existence of Auslander-Reiten sequences in \(\mod_A\).

\begin{thm}\label{THM: existence of AR sequences in mod_A for all A}
     Let \(M\) be an indecomposable finitely generated non-projective \(A\)-module and let \(L\) an indecomposable finitely generated non-injective \(A\)-module. Let \(\operatorname{t}\) and \(\operatorname{t}^-\) denote the standard Auslander-Reiten translates in \(\mod_A\) and let \(\operatorname{N}^-(-)=\Hom_A(\D(A),-)\) be the standard Nakayama functor on \(\mod_A\). Then the following hold:
    \begin{itemize}
        \item [1.] There is an Auslander-Reiten sequence in \(\mod_A\) of the form
        \[
        0 \xrightarrow{} \operatorname{t}(M) \xrightarrow{} E_M \xrightarrow{}M \xrightarrow{} 0.
        \]
        \item [2.] There is an Auslander-Reiten sequence in \(\mod_A\) of the form
        \[
        0 \xrightarrow{} L \xrightarrow{} F_L \xrightarrow{} \operatorname{t}^-(L) \xrightarrow{} 0.
        \]
    \end{itemize}
\begin{proof}
    \textit{Part 1} and \textit{Part 2.} By \cite[Theorem on page~16 and Lemma on page~62]{happel-triangulated-1988}, the category \(\underline{\mod}_{\widehat{A}}\) is a triangulated category. We denote the suspension functor of \(\underline{\mod}_{\widehat{A}}\) by \(\Sigma\). By \cite[Proposition on page~67]{happel-triangulated-1988}, \(\mod_A\) is the heart of a t-structure on \(\underline{\mod}_{\widehat{A}}\) and hence by \cite[Theorem~1.3.6 on page~31]{beilinson-faisceaux-1982}, \(\mod_A\) is an extension closed \(k\)-linear proper abelian subcategory of \(\underline{\mod}_{\widehat{A}}\). By \cite[Lemma on page~63]{happel-triangulated-1988}, we have isomorphisms \(\underline{\Hom}_{\widehat{A}}(M,\Sigma^i L) \cong \Ext^i_A(M,L)\) for all \(i\) in \(\bbZ\) and therefore, \(\underline{\Hom}_{\widehat{A}}(\mod_A,\Sigma^{-1}\mod_A)=0\). It was known in \cite[Lemma on page~354]{hughes-trivial-1983}, that \(\mod_{\widehat{A}}\) has Auslander-Reiten sequences. Later, it was stated in \cite{happel-triangulated-1988}, that \(\underline{\mod}_{\widehat{A}}\) has the so-called Auslander-Reiten triangles (one can use \cite[Proposition~5.11]{iyama-auslander-2024}, for example). Then by \cite[Theorem~I.2.4]{reiten-noetherian-2002} (see also Remark~\ref{REM: removing algebraically closed field}), the triangulated category \(\underline{\mod}_{\widehat{A}}\) has a Serre functor. As \(\mod_{\widehat{A}}\) is Hom-finite, so is \(\underline{\mod}_{\widehat{A}}\). To show that \(\underline{\mod}_{\widehat{A}}\) is Krull-Schmidt, it suffices to show that each idempotent in \(\underline{\mod}_{\widehat{A}}\) splits (see \cite[page~52]{ringel-tame-1984}, \cite[page~26]{happel-triangulated-1988}, \cite[Corollary~A.2]{chen-algebras-2008}, \cite[Section~4]{krause-krull-2015} and \cite[Theorem~6.1]{shah-krull-2023}). 
    
    Let us first show that each idempotent in \(\mod_{\widehat{A}}\) splits. Following \cite[page~60]{happel-triangulated-1988}, we let \(M=(M_i,m_i)_{i\in\bbZ}\) be a finitely generated \(\widehat{A}\)-module. That is, \(M_i\) is a finitely generated \(A\)-module with finitely many \(M_i\) being nonzero and \(M_i \xrightarrow{m_i} \N^-(M_{i+1})\) is an \(A\)-module homomorphism such that \(\N^-(m_{i+1}) m_i=0\). Let \(M \xrightarrow{e=(e_i)_{i\in\bbZ}} M\) be an \(\widehat{A}\)-module endomorphism. That is, \(M_i \xrightarrow{e_i} M_i\) is an \(A\)-module homomorphism such that
    \[
    \begin{tikzcd}
        M_i \arrow[r, "m_i"] \arrow[d, "e_i"'] & \N^-(M_{i+1}) \arrow[d, "\N^-(e_{i+1})"] \\
        M_i \arrow[r, "m_i"]                   & \N^-(M_{i+1})                             
    \end{tikzcd}
    \]
    commutes for each \(i\) in \(\bbZ\). Assume that \(e\) is an idempotent in \(\End_{\widehat{A}}(M)\). Then for each \(i\) in \(\bbZ\), the \(A\)-module homomorphism \(M_i \xrightarrow{e_i} M_i\) is an idempotent in \(\End_A(M_i)\). Hence, for each \(i\) in \(\bbZ\), there exists a factorisation \(e_i\colon M_i \xrightarrow{r_i} X_i \xrightarrow{s_i} M_i\) of \(e_i\) such that \(r_i\) and \(s_i\) are \(A\)-module homomorphisms satisfying the equation \(1_{X_i}=r_is_i\). We define an \(A\)-module homomorphism \(X_i \xrightarrow{x_i} \N^{-}(X_{i+1})\) as the following composition:
    \[
    x_i\colon X_i \xrightarrow{s_i} M_i \xrightarrow{m_i} \N^{-}(M_{i+1}) \xrightarrow{\N^{-}(r_{i+1})} \N^{-}(X_{i+1}).
    \]
    It is not hard to check that the pair \((X_i,x_i)_{i\in\bbZ}\) defines a finitely generated \(\widehat{A}\)-module and that the sequences \(r=(r_i)_{i\in\bbZ}\) and \(s=(s_i)_{i\in\bbZ}\) define \(\widehat{A}\)-module homomorphisms such that \(e=sr\) and \(1_M=rs\). Hence, \(e\) splits. By \cite[Proposition~5.9]{cortes-izurdiaga-reflective-2023}, the stable module category \(\underline{\mod}_{\widehat{A}}\) has split idempotents and therefore is Krull-Schmidt. Both parts of the theorem now follow from Proposition~\ref{THM: AR sequences in modA}.
\end{proof}
\end{thm}

\paragraph{\bf Acknowledgement.} I am grateful to Peter Jørgensen for his invaluable supervision. I deeply appreciate our insightful discussions, which have significantly enhanced my understanding of algebra. I am also immensely thankful to everyone at the (extended) Aarhus Homological Algebra Group. Special thanks to David Pauksztello for his useful feedback on the first draft of this article and to the kind referees for their comments and suggestions. Everyone mentioned above has helped improve the readability of this article and deepened my understanding of the topic. This work was supported by the Independent Research Fund Denmark (grant no. 1026-00050B).

\footnotesize
\bibliographystyle{alpha}
\bibliography{references}

\textsc{Department of Mathematics, Aarhus University, Ny Munkegade 118, 8000 Aarhus C, Denmark}.\\
Email address: \href{mailto:david.nkansah@math.au.dk}{david.nkansah@math.au.dk}

\end{document}